\documentstyle[12pt]{amsart}
\begin{document}

\title{A combinatorial approach to quantification of Lie algebras}

\author{V.K. Kharchenko}
\address{\hskip-\parindent
        V.K. Kharchenko\\
        UNAM\\
        Calle de 1 mayo s/n, 54769\\
        Edo. de M\'exico, M\'exico.\\
        and Sobolev Institute of Mathematics\\
        Novosibirsk, 630090 \\
        Russia}
\email{vlad\@ servidor.unam.mx}

\thanks{Research at MSRI is supported in part by NSF grant DMS-9701755.
Supported in part by CONACyT M\'exico, grant 32130-E}

\begin{abstract}
We propose a notion of a quantum universal enveloping algebra for an
arbitrary Lie algebra defined by generators and relations which is
based
on the quantum Lie  operation concept. This enveloping algebra has a
PBW
basis that admits the Kashiwara crystalization. We describe all skew
primitive elements of the quantum universal enveloping algebra for the
classical  nilpotent algebras of the infinite series defined by the
Serre
relations  and prove that the set of PBW-generators for each of these
enveloping algebras coincides with the Lalonde--Ram basis of the
ground
Lie algebra  with a skew commutator in place of the Lie operation. The
similar statement is valid for Hall--Shirshov basis of any Lie algebra
defined  by one relation, but it is not so in general case.
\end{abstract}

\maketitle
\newcounter{par}
\newcounter{nom}[par]
\newtheorem{nom}[par]{Nom}
\newtheorem{theorem}{\sc Theorem}
\newtheorem{corollary}{\sc Corollary}
\newtheorem{lemma}{\sc Lemma}
\newtheorem{definition}{\sc Definition}
\renewcommand{\thetheorem}{\thepar.\arabic{nom}}
\renewcommand{\thelemma}{\thepar.\arabic{nom}}
\renewcommand{\thecorollary}{\thepar.\arabic{nom}}
\renewcommand{\thedefinition}{\thepar.\arabic{nom}}

\centerline{1.\large  \it Introduction}

\

\vspace{-1mm}
\setcounter{par}{1}
\noindent Quantum universal enveloping algebras appeared in the famous papers by Drinfeld
[{\bf \ref{pDri}}] and Jimbo [{\bf \ref{pJim}}]. Since then a great deal of articles and number
of monographs were devoted to their investigation. All of these researches are mainly
concerned with  a particular quantification of Lie algebras of the classical series. This is
accounted for first by the fact that these Lie algebras have applications  and visual
interpretations in physical speculations, and then by the fact that a general, and
commonly accepted as standard, notion of a quantum universal enveloping algebra is not
elaborated yet (see a detailed discussion in [{\bf \ref{pBau}},
{\bf \ref{pLyu}}]).

In the present paper we propose a combinatorial solution of this problem by means of the
quantum (Lie) operation concept [{\bf \ref{pKh}, \ref{pJA}, \ref{pDAN}}].
In line with the main idea of
our approach,
the skew primitive elements must play the same role
in quantum enveloping algebras  as
the primitive elements do in the classical case. By the Friedrichs
criteria [{\bf \ref{pCohn}}, {\bf \ref{pFri}},
{\bf \ref{pLyn}}, {\bf \ref{pMag}}, {\bf \ref{pMMo}}], the primitive
elements form the ground Lie algebra in the classical case. For this reason we consider
the space spanned by the skew primitive elements and equipped with the quantum
operations as a quantum analogue of a Lie algebra.

In the second section we adduce the main notions and consider some examples. These
examples, in particular,  show that the Drinfeld--Jimbo enveloping algebra as well as its
modifications are quantum enveloping algebras in our sense.

In the third section with the help of the Heyneman--Radford theorem we
introduce a notion of a {\it combinatorial rank} of a Hopf algebra
generated by
skew primitive semi-invariants. Then we define the quantum enveloping algebra of an
{\it arbitrary rank} that slightly generalises  the definitions given in  the preceding
section.

The basis construction problem for the quantum enveloping algebras is considered in the
fourth section. We indicate two main methods for the construction of {\it PBW-generators}.
One of them modifies  the Hall--Shirshov  basis construction process by
means of  replacing  the Lie operation with a skew commutator. The
set of the PBW-generators defined in this way,
the values of {\it hard super-letters}, plays the
same role as the basis of the ground Lie algebra does in the PBW theorem. At first
glance it would seem reasonable to consider the ${\bf k}[G] $-module generated by
the values of hard super-letters as a quantum Lie algebra. However, this extremely
important module falls far short of  being uniquely defined. It essentially depends on the
ordering of the main generators, their degrees, and it is almost never antipode stable.
Also we have to note the following important fact.
Our definition of the hard super-letter is not constructive and, of course, it cannot be constructive
in general. The basis construction problem includs the word problem for Lie algebras
defined by generators and relations, while the latter one has no general algorithmic solution
(see [{\bf \ref{pBo}}, {\bf \ref{pBK}}]).

The second method is connected with the Kashiwara
crystallisation idea [{\bf \ref{pKas}}, {\bf \ref{pKas1}}]
(see also a development
in [{\bf \ref{pCli}}, {\bf \ref{pKit}}]). M. Kashiwara has considered the
main parameter $q$ of the Drinfeld--Jimbo  enveloping
algebra as a temperature of some physical medium. When the temperature tend to zero,
the medium crystallises. The PBW-generators must crystallise  as well.
In our case under this process   no one limit quantum enveloping algebra appears since
the  existence conditions normally include equalities of the form $\prod p_{ij}=1$
(see [{\bf \ref{pJA}}]).
Nevertheless if we equate all quantification parameters to zero, the hard super-letters
would form a new set of PBW-generators for the given quantum universal
enveloping algebra. To put this another way, the PBW-basis defined by the
super-letters admits the Kashiwara crystallisation.

In the fifth section we bring a way to construct a Groebner--Shirshov   relations system
for a quantum enveloping algebra. This system is related to the main skew
primitive generators, and, according to the Diamond Lemma
(see [{\bf \ref{pBer}}, {\bf \ref{pBok}}, {\bf \ref{pSh2}}]),
it determines the crystal basis.
The usefulness of the
Groebner--Shirshov systems depends upon the fact that such
a system not only defines a
basis of an associative algebra, but it also provides a
simple diminishing algorithm for
expansion of elements on this basis (see, for example [{\bf \ref{pBMM}}]).

In the sixth section we adapt  a well known method of triangular
splitting to the
{\it quantification with constants}.
The original method  appeared in studies of simple
finite dimensional  Lie algebras. Then it has been extended into
the field of quantum algebra in a lot of publications
(see, for example [{\bf \ref{pBM}},
{\bf \ref{pLuz}}, {\bf \ref{pYa}}]).
By means of this method the investigation of the
Drinfeld--Jimbo  enveloping algebra amounts to
a consideration of its positive and
negative homogeneous components,
{\it quantum Borel sub-algebras}.

In the seventh section we consider  more thoroughly
the quantum universal
enveloping algebras of nilpotent
algebras of the series $A_n,$ $B_n,$ $C_n,$ $D_n$
defined by the Serre relations.
We adduce first lists of all
hard super-letters in the explicit form, then
Groebner--Shirshov relations systems, and next
spaces $L(U_P(\frak g))$ spanned by the
skew primitive elements (i.e. the Lie algebra
quantifications ${\frak g}_P$
proper). In all cases the lists of hard super-letters
(but the hard super-letters themselves)
turn out to be independent of the quantification parameters.
This means that the
PBW-generators result from  the Hall--Shirshov
basis of the ground Lie algebra by
replacing the Lie operation with the skew commutator.
The same is valid for the
Groebner--Shirshov  relations systems. Note that the
Hall--Shirshov bases, under the name {\it standard Lyndon bases},
for the classical Lie series were
constructed by P. Lalonde and A. Ram [{\bf \ref{pLR}}], while the
Groebner--Shirshov systems of Lie
relations were
found by L.A. Bokut' and A.A. Klein  [{\bf \ref{pKle}}].

Furthermore, in all cases ${\frak g}_P$
as a quantum Lie algebra (in our sense) proves to be very simple in structure. Either it is
a coloured  Lie super-algebra (provided that the main parameter $p_{11}$ equals
1), or values of all non-unary  quantum operations equal zero on ${\frak g}_P.$ In
particular, if char({\bf k}) $=0$ and $ p_{11}^t\neq 1$ then the partial quantum
operations may be defined on  ${\frak g}_P,$ but all of them have zero values.
Thus, in this case we have a reason to consider $U_P({\frak g})$ as an algebra of
`commutative'  quantum polynomials, since the universal enveloping algebra of a Lie
algebra with zero bracket is the algebra of ordinary commutative polynomials.
From this standpoint the Drinfeld--Jimbo enveloping algebra is a  `quantum' Weyl algebra
of (skew) differential operators. Immediately afterwards a number of interesting
questions appears.
What is the structure of other algebras of  `commutative' quantum
polynomials? Under what conditions are the
quantum universal enveloping algebras of
homogeneous components of other Kac--Moody  algebras
defined by the Gabber--Kac
relations [{\bf \ref{pKac}}] the algebras of  `commutative'  quantum
polynomials? When do the PBW-generators result
from a basis of the ground Lie algebra by means of replacing the Lie
operation with the skew commutator? These and
other questions we briefly discuss  in the last section.

It is as well to bear in mind that the combinatorial
approach is not free from flaws: the quantum universal enveloping
algebra essentially depends on a
combinatorial representation of the ground Lie algebra,
i.e. a close connection with the
abstract category of Lie algebras is lost.

\

\centerline{2. \large \it Quantum enveloping algebras}

\

\vspace{-1mm}
\stepcounter{par}
\noindent Recall that a variable $x$
is called a {quantum variable} if an element $g_x$
of a fixed  Abelian group $G$
and a character
$\chi ^x\in G^*$
are associated with it.
A noncommutative polynomial in quantum variables is called
a {\it quantum operation} if all of its values in all Hopf algebras are
skew primitive provided that every variable $x$
has a value $x=a$
such that
\begin{equation}
\Delta (a)=a\otimes 1+g_x\otimes a,\ \ \
g^{-1}ag=\chi ^x(g)a,  \ \ \ g\in G.
\label{AI}
\end{equation}

Let $x_1,\ldots ,x_n$
be a set of quantum variables. For each word $u$
in $x_1,\ldots ,x_n$
we denote by $g_u$
an element of $G$
that appears from $u$
by replacing  of all $x_i$ with $g_{x_i}.$
In the same way we denote by $\chi ^u$
a character that appears from $u$
by replacing of all $x_i$
with $\chi ^{x_i}.$
Thus on the free algebra
{\bf k}$\langle x_1,\ldots , x_n\rangle$
a grading by the group $G\times G^*$
is defined.
For each pair of homogeneous elements $u,$\ $v$
we fix the denotations
$p_{uv}=\chi ^u(g_v)=p(u,v).$

The quantum operation can be defined equivalently as a $G\times 1$-homogeneous
polynomial  that has only primitive values
in all {\it braided bigraded Hopf algebras} provided that all quantum
variables have primitive homogeneous values
$g_a=g_x,$ $\chi ^a=\chi ^x$
(see [{\bf \ref{pKh}}, Sect. 1--4]).

Recall that a {\it constitution} of a word $u$
is a sequence of non-negative integers $(m_1,m_2,\ldots , m_n)$
such that $u$
is of degree $m_1$ in $x_1,$ deg$_1(u)=m_1;$
of degree $m_2$ in $x_2,$ deg$_2(u)=m_2;$
and so on. Since the group $G$
is abelian, all constitution homogeneous polynomials are
homogeneous with respect to the grading. Let us
define a bilinear skew commutator on the set of
graded homogeneous noncommutative  polynomials by the formula
\begin{equation}
[u,v]=uv-p_{uv}vu.
\label{op1}
\end{equation}
These brackets satisfy the following Jacobi
and skew  differential identities:
\begin{equation}
[[u,v],w]=[u,[v,w]]+p_{wv}^{-1}[[u,w],v]+(p_{vw}-p_{wv}^{-1})[u,w]\cdot v;
\label{Jak1}
\end{equation}
\begin{equation}
[[u,v],w]=[u,[v,w]]+p_{vw}[[u,w],v]+p_{uv}(p_{vw}p_{wv}-1)v\cdot [u,w];
\label{Jak2}
\end{equation}
\begin{equation}
[u,v\cdot w]=[u,v]\cdot w+p_{uv}v\cdot [u,w];  \  \
[u\cdot v,w]=p_{vw}[u,w]\cdot v+u\cdot [v,w],
\label{dif1}
\end{equation}
where by the dot we denote the usual multiplication.
It is easy to see that the following conditional
{\it restricted identities} are valid as well
\begin{equation}
[u,v^n]=[\ldots [[u,v],v]\ldots ,v];  \  \  \
[v^n,u]=[v,[\ldots [v,u]\ldots ]],
\label{ogr}
\end{equation}
provided that $p_{vv}$
is a primitive $t$-th
root of unit, and $n=t$
or $n=tl^k$ in the case of characteristic $l>0.$

Suppose that a Lie algebra $\frak g$
is defined by the generators $x_1,\ldots ,x_n$
and the relations $f_i=0.$
Let us convert the generators into quantum variables. For this associate
to them elements of  $G\times G^*$
in arbitrary way. Let $P=||p_{ij}||,$  $p_{ij}=\chi ^{x_i}(g_{x_j})$
be the {\it quantification matrix}.
\addtocounter{nom}{1}
\begin{definition} \rm
A {\it braided quantum enveloping algebra} is a braided
bigraded Hopf algebra
$U_P^b(\frak g)$
defined by the variables
$x_1,\ldots ,x_n$
and the relations $f_i=0,$
where the Lie operation is replaced with (\ref{op1}),
provided that in this way $f_i$
are converted into the quantum operations $f_i^*.$
The coproduct and the commutation relations in the tensor product are defined by
\begin{equation}
\Delta (x_i)=x_i\underline{\otimes }1+1\underline{\otimes } x_i,
\label{bcop}
\end{equation}
\begin{equation}
(x_i\underline{\otimes }x_j)\cdot (x_k\underline{\otimes }x_m)=
(\chi ^{x_k}(g_{x_j}))^{-1}x_ix_k\underline{\otimes }x_jx_m.
\label{bcom}
\end{equation}
\label{qvabr}
\end{definition}
\addtocounter{nom}{1}
\begin{definition} \rm
A {\it simple quantification} or a
{\it quantum universal enveloping algebra} of
$\frak g$
is an algebra $U_P({\frak g})$
that is isomorphic to the skew group algebra
\begin{equation}
U_P({\frak g})=U_P^b({\frak g})*G,
\label{ra1}
\end{equation}
where the group action and the coproduct are defined by
\begin{equation}
g^{-1}x_ig=\chi ^{x_i}(g)x_i,\
\Delta (x_i)=x_i\otimes 1+g_{x_i}\otimes x_i, \ \Delta (g)=g\otimes g.
\label{ccop}
\end{equation}
\label{qva}
\end{definition}
\addtocounter{nom}{1}
\begin{definition} \rm
A {\it quantification with constants} is a simple quantification
where additionally some generators $x_i$
associated to the trivial character are replaced with
the constants $\alpha _i(1-g_{x_i}).$
\label{qvaco}
\end{definition}

The formulae (\ref{ccop}) and (\ref{bcop}) correctly define the
coproduct since by definition of the quantum operation
$\Delta(f_i^*)=f_i^*\otimes 1+g\otimes f_i^*$
in the case of ordinary Hopf algebras and
$\Delta(f_i^*)=f_i^*\underline{\otimes }1+1\underline{\otimes }f_i^*$
in the braided case.

We have to note that the defined quantifications essentially depend on
the combinatorial representation of the Lie algebra. For example,
an additional relation $[x_1,x_1]=0$
does not change the Lie algebra. At the same time if
$\chi ^{x_1}(g_1)=-1$
then this relation admits the quantification and yields a nontrivial
relation for the quantum enveloping algebra, $2x^2_1=0.$

{\sc Example 1.} Suppose that the Lie algebra is defined by a system
of constitution homogeneous relations. If the characters $\chi ^i$
are such that $p_{ij}p_{ji}=1$
for all $i,j$
then the skew commutator itself is a quantum operation.
Therefore on replacing the Lie operation
all relations become quantum operations as well.
This means that the braided enveloping algebra is the universal
enveloping algebra $U({\frak g}^{col})$
of the coloured Lie super-algebra which is defined by the same relations
as the given Lie algebra is.
The simple quantification appears as the Radford biproduct
$U({\frak g}^{col})\star {\bf k}[G]$
or, equivalently, as the universal $G$-enveloping algebra
of the coloured Lie super-algebra ${\frak g}^{col}$
(see [{\bf \ref{pRad}}] or [{\bf \ref{pKh}}, Example 1.9]).

{\sc Example 2.} Suppose that the Lie algebra  $\frak{g}$
is defined by the generators  $x_1,\ldots , x_n$
and the system of nil relations
\begin{equation}
x_j(ad x_i)^{n_{ij}}=0, \ \ \ 1\leq i\neq j\leq n.
\label{rel}
\end{equation}

Usually instead of the matrix of degrees (without the main diagonal)
$||n_{ij}||$
the matrix $A=||a_{ij}||,$ $a_{ij}=1-n_{ij}$
is considered. The Coxeter graph ${\bf \Gamma }(A)$
is associated to every such a matrix.
This graph has the vertices $1,\ldots ,n,$
where the vertex  $i$
is connected by $a_{ij}a_{ji}$
edges with the vertex $j.$

If $a_{ij}=0$
then the relation $x_j$ad$x_i=0$
is in the list (\ref{rel}), and the relation $x_i$(ad$x_j)^{n_{ji}}=0$
is a consequence of it. The skew commutator $[x_j,x_i]$
is a quantum operation if and only if $p_{ij}p_{ji}=1.$
Under this condition we have $[x_i,x_j]=-p_{ij}[x_j,x_i].$
Therefore both in the given Lie algebra and in its
quantification one may replace the relation $x_i$(ad$x_j)^{n_{ji}}=0$
with $x_i$ad$x_j=0.$
In other words, without loss of generality, we may suppose that
$a_{ij}=0 \leftrightarrow a_{ji}=0.$
By the Gabber-Kac theorem [{\bf \ref{pKac}}]
we get that the algebra $\frak g$
is the positive homogeneous component
${\frak{g}}_1^+$
of a Kac-Moody algebra ${\frak{g}}_1.$

Theorem 6.1 [{\bf \ref{pKh}}] describes the conditions  for a homogeneous
polynomial in two variables which is linear in one of them to be a quantum operation.
From this theorem we have the following corollary.
\addtocounter{nom}{1}
\begin{corollary}
If $n_{ij}$
is a simple number or unit and in the former case  $p_{ii}$
is not a primitive $n_{ij}$-th root of unit,
then the relation $(\ref{rel})$
admits a quantification if and only if $p_{ij}p_{ji}=p^{a_{ij}}_{ii}.$
\label{sled}
\end{corollary}
Theorem 6.1 [{\bf \ref{pKh}}] provides no essential restrictions on the non-diagonal
parameters $p_{ij}:$
if the matrix $P$
correctly defines a quantification of (\ref{rel}) then for every set
$Z=\{ z_{ij}|z_{ij}z_{ji}=z_{ii}=1\} $
the following matrix does as well:
 \begin{equation}
P_Z=\{ p_{ij}z_{ij}|p_{ij}\in P,z_{ij}\in Z\}.
\label{col}
\end{equation}

{\sc Example 3}. Let $G$
be freely generated by $g_1,\ldots g_n$
and $A$
be a generalised Cartan matrix symmetrised by  $d_1,\ldots ,d_n,$
while the characters are defined by $p_{ij}=q^{-d_ia_{ij}}.$
In this case the simple quantification is the positive component
of the Drinfeld--Jimbo enveloping algebra  together with
the group-like elements, $U_P({\frak g})= U_q^+({\frak{g}})*G.$
By means of an arbitrary  deformation (\ref{col}) one may
define a `colouring' of  $ U_q^+({\frak{g}})*G.$

The braided enveloping algebra equals  $U_q^+({\frak{g}})$
where the coproduct and braiding are defined by (\ref{bcop})
and (\ref{bcom}) with the coefficient $q^{d_ka_{kj}}.$
The formula (\ref{col}) correctly defines  its `colouring' as well.

{\sc Example 4}. If in the example above we complete the set
of quantum variables by the new ones
$x_1^-,\ldots ,x_n^-;$ $z_1,\ldots ,z_n$
such that
\begin{equation}
\chi ^{x^-}=(\chi ^x)^{-1},\ \ g_{x^-}=g_x,\ \
\chi ^{z}=\hbox{id},\ \ g_{z_i}=g_i^2,
\label{skon}
\end{equation}
then by  [{\bf \ref{pKh}}, Theorem 6.1]
the Gabber--Kac relations (2), (3), and $[e_i,f_j]=\delta _{ij}h_i$
(see [{\bf \ref{pKac}}, Theorem 2]) under  the identification
$e_i=x_i,$ $f_i=x_i^-,$ $h_i=z_i$
admit the quantification with constants
$z_i=\varepsilon _i\, (1-g_i^2).$ (Unformally we may consider
the obtained quantification as one of the Kac--Moody algebra identifying
$g_i$ with  $q^{h_i},$
where
the rest of the Kac--Moody algebra relations, $[h_i,e_j]=a_{ij}e_i,$ $[h_i, f_j]=-a_{ij}f_j,$
is quantified to the $G$-action:  $g_j^{-1}x_i^{\pm }g_j=q^{\mp d_{ij}a_{ij}}x_i^{\pm }.$)
This quantification coincides with
the Drinfeld--Jimbo  one
under a suitable choice of $x_i,$ $x_i^-,$ and $\varepsilon _i$
depending up the particular definition of $U_q({\frak g}):$
$$
\begin{matrix}
[{\bf \ref{pLu}}] \hfill &x_i=E_i,\ g_i=K_i,\ x_i^-=F_iK_i, \
p_{ij}=v^{-d_ia_{ij}}, \
\varepsilon _i=(v^{-d_i}-v^{d_i})^{-1}; \hfill \cr
[{\bf \ref{pLuz}}]
 \hfill &x_i=E_i,\ g_i=\tilde{K}_i,
\ x_i^-=F_i\tilde{K}_i, \
p_{i\mu }=v^{-\langle \mu ,i^{\prime }\rangle }, \
\varepsilon _i=(v_i^{-1}-v_i)^{-1}; \hfill \cr
[{\bf \ref{pKas}}] \Delta _+\hfill &x_i=e_i,\ g_i=t_i,\ x_i^-=t_if_i, \
p_{ij}=q_j^{-\langle h_j,\alpha _i\rangle },\
\varepsilon _i=(q_i-q_i^3)^{-1};\hfill \cr
[{\bf \ref{pKas}}] \Delta _-\hfill &x_i=f_i,\ g_i=t_i,\ x_i^-=e_it_i, \
p_{ij}=q_j^{\langle h_j,\alpha _i\rangle },\
\varepsilon _i=(q_i^{-1}-q_i)^{-1};\hfill \cr
[{\bf \ref{pMon}}] \hfill &x_i=E_iK_i, \ g_i=K_i^2, \ x_i^-=F_iK_i,\
p_{ij}=q^{-2d_ia_{ij}}, \ \varepsilon _i=(1-q^{4d_i})^{-1}.\hfill
\end{matrix}
$$
By (\ref{skon}) the brackets $[x_i,x_j^-]$
are  quantum operations only if $p_{ij}=p_{ji}.$
So in this case the `colourings' (\ref{col}) may
be only black-white, $z_{ij}=\pm 1.$

In the perfect analogy the Kang quantification [{\bf \ref{pKang}}] of the generalised Kac-Moody algebras
[{\bf \ref{pBor}}] is a quantification  in our sense as well.

\

\centerline{3.\large \it Combinatorial rank}

\

\vspace{-1mm}
\stepcounter{par}
\addtocounter{nom}{1}
\noindent  By the above definitions the quantum enveloping algebras
(with or without constants) are character Hopf algebras
(see [{\bf \ref{pKh}}, Definition 1.2]). In this section
by means of a  combinatorial rank notion we
identify the quantum enveloping algebras
in the class of  character Hopf algebras.

Let $H$
be a character Hopf algebra generated by $a_1,\ldots ,a_n:$
\begin{equation}
\Delta (a_i)=a_i\otimes 1+g_{a_i}\otimes a_i,\ \ \ \
g^{-1}a_ig=\chi ^{a_i}(g)a_i, \ \ g\in G.
\label{AIc}
\end{equation}
Let us associate a quantum variable $x_i$
with the parameters $(\chi ^{a_i},$ $g_{a_i})$
to $a_i.$
Denote by $G\langle X\rangle$
the free enveloping algebra defined by the quantum variables
$x_1\ldots ,x_n.$
(see [{\bf \ref{pKh}}, Sec. 3] under denotation $H\langle X\rangle $).

The map $x_i\rightarrow a_i$
has an extension to a homomorphism of Hopf algebras
$\varphi : G\langle X\rangle \rightarrow H.$
Denote by $I$
the kernel of this homomorphism. If $I\neq 0$
then by the Heyneman--Radford theorem  (see, for example
[{\bf \ref{pMon}}, pages  65--71]),
the Hopf ideal $I$
has a non-zero skew primitive element. Let $I_1$
be an ideal generated by all skew primitive elements of $I.$
Clearly $I_1$
is a Hopf ideal as well. Now consider the Hopf ideal $I/I_1$
of the quotient Hopf algebra $G\langle X\rangle /I_1.$
This ideal also has non-zero skew primitive elements (provided $I_1\neq I$).
Denote by $I_2/I_1$
the ideal generated by all skew primitive elements of $I/I_1,$
where $I_2$ is
its preimage with respect to the projection
$G\langle X\rangle \rightarrow G\langle X\rangle /I_1.$
Continuing the process we will find a strictly increasing,
finite or infinite, chain of Hopf ideals of $G\langle X\rangle :$
\begin{equation}
0=I_0\subset I_1\subset I_2 \subset \ldots \subset I_n\subset \ldots .
\label{tce}
\end{equation}
\begin{definition}\rm The length
of (\ref{tce}) is called a {\it combinatorial rank} of $H.$
\label{ra}
\end{definition}
By definition, the combinatorial rank of any quantum
enveloping algebra (with constants) equals one. In the case of zero
characteristic the inverse statement is valid as well.
\addtocounter{nom}{1}
\begin{theorem}
Each character Hopf algebra of the combinatorial rank 1
over a field of zero characteristic is isomorphic to a
quantum enveloping algebra with constants of a Lie algebra.
\label{kh}
\end{theorem}
{\it Proof}. By definition, $I$
is generated by skew primitive
elements. These elements as noncommutative polynomials are the quantum operations.
Consider one of  them, say $f.$
Let us decompose $f$
into a sum of homogeneous components $f=\sum f_i.$
All positive components belongs to ${\bf k}\langle X\rangle $
and they are the quantum operations themselves,
while the constant component has the form $\alpha (1-g), g\in G$
(see [{\bf \ref{pKh}}, Sec. 3 and Prop. 3.3]). If $\alpha \neq 0$
then we introduce a new quantum variable $z_f$
with the parameters $(id, g).$
Each $f_i$
has a representation  through the skew commutator.
Indeed, by [{\bf \ref{pKh}}, Theorem 7.5]
the complete linearization $f_i^{lin}$
of $f_i$
has the required representation. By the identification of
variables in a suitable way in $f_i^{lin}$
we get the required representation for $f_i$
multiplied by a natural number, $m_if_i=f_i^{[\,]}.$

Now consider a Lie algebra ${\frak g}$
defined by the generators $x_i, z_f$
and the relations $\sum m_i^{-1}f_i^{[\,]}+z_f=0,$
with the Lie multiplication in place of
the skew commutator.
It is clear that $H$
is the quantification with constants of ${\frak g}.$

In the same way one may introduce the notion of the combinatorial rank
for the braided bigraded Hopf algebras.  In this case
all braided quantum enveloping algebras are of rank 1,
and all braided bigraded algebras of rank 1 are
the braided quantification of some Lie algebras.

Now we are ready to define a quantification of {\it arbitrary rank}.
For this in the definitions of the above section it is necessary
to change the requirement that all $f_i^*$
are quantum operations with the following  condition.

{\it The set $F$
splits in a union $F=\cup_{j=1}^nF_j$
such that  $F_1^*$
consists of  quantum operations; the set $F_2^*$
consists of  skew primitive elements of $G\langle X||F_1^*\rangle ;$
the set $F_3^*$
consists of skew primitive elements of $G\langle X||F_1^*, F_2^*\rangle ,$
and so on}.

The quantum enveloping algebras of an arbitrary
rank are character Hopf algebras also.
But it is not clear if any character Hopf algebra
is a quantification of some rank of a suitable Lie algebra.
It is so if the Hopf algebra is homogeneous and the ground
field has a zero characteristic (to appear).
Also it is not clear if there exist  character Hopf algebras,
or braided bigraded Hopf algebras, of infinite combinatorial rank;
while it is easy to see that $\cup I_n=I.$
Also it is possible to show that $F_1$
always contains all relations of a minimal constitution in $F.$
For example, each of (\ref{rel})
is of a minimal constitution in (\ref{rel}).
Therefore the quantification of arbitrary rank with the identification $g_i=\hbox{exp}(h_i)$
of any (generalized) Kac--Moody algebra ${\frak g},$
or its nilpotent component ${\frak g}^+$,
is always a quantification in the sense of the above section.

\

\centerline{4.\large  \it PBW-generators and crystallisation}

\

\vspace{-1mm}
\stepcounter{par}
\noindent
The next result yields a PBW basis
for  the  quantum enveloping algebras.
\addtocounter{nom}{1}
\begin{theorem}
Every character Hopf algebra $H$
has a linearly ordered  set of constitution homogeneous elements
$U=\{ u_i\ |\ i\in I\}$
such that the set of all products $gu_1^{n_1}u_2^{n_2}\cdots u_m^{n_m},$
where $g\in G,$ $ u_1<u_2<\ldots <u_m,$ $0\leq n_i<h(i)$
forms a basis of $H.$
Here if $p_{ii}\stackrel{df}{=}p_{u_i u_i}$
is not a root of unity then $h(i)=\infty ;$
if $ p_{ii}=1$
then either $h(i)=\infty $
or $h(i)=l$
is the characteristic of the ground field;
if $ p_{ii}$
is a primitive $t$-th
root of unity, $t\neq 1,$
then $h(i)=t.$
\label{can}
\end{theorem}
The set $U$
is referred to as  a set of {\it PBW-generators} of $H.$
This theorem  easily follows from [{\bf \ref{pKh1}}, Theorem 2].
Let us recall necessary notions.

Let $a_1,\ldots ,a_n$
be a set of skew primitive generators of $H,$
and let $x_i$
be the associated quantum variables.
Consider the lexicographical ordering of all
words in $x_1>x_2>\ldots >x_n.$
A non-empty word $u$
is called {\it standard} if $vw>wv$
for each decomposition $u=vw$
with non-empty $v,w.$
The following properties are well known
(see, for example  [{\bf \ref{pFox}}, {\bf \ref{pCoh}},
{\bf \ref{pLot}}, {\bf \ref{pSh1}}, {\bf \ref{pSh2}}]).

1s. A word  $u$
is standard if and only if it is greater
than each of its ends.

2s. Every standard word starts with a maximal letter that it has.

3s. Each word $c$
has a unique representation  $c=u_1^{n_1}u_2^{n_2}\cdots u_k^{n_k},$
where $u_1<u_2<\cdots <u_k$
are standard words (the Lyndon theorem).

4s. If $u,v$
are different standard words and $u^n$
contains $v^k$
as a sub-word, $u^n=cv^kd,$
then $u$
itself contains $v^k$
as a sub-word, $u=bv^ke.$

The set of {\it standard nonassociative} words is defined
as the smallest set $SL$
that contains all variables $x_i$
and satisfies the following properties.

1)\ If $[u]=[[v][w]]\in SL$
then $[v],[w]\in SL,$
and $v>w$
are standard.

2)\  If $[u]=[[[v_1][v_2]][w]]\in SL$ then $v_2\leq w.$

The following statements are valid as well.

5s. Every standard word has
the only alignment of brackets such that the appeared
nonassociative word is standard (the Shirshov theorem [{\bf \ref{pSh1}}]).

6s. The factors $v, w$
of the nonassociative decomposition $[u]=[[v][w]]$
are the standard words such that $u=vw$
and  $v$
has the minimal length ([{\bf \ref{pSh2}}]).
\addtocounter{nom}{1}
\begin{definition} \rm A {\it super-letter}
is a polynomial that equals a nonassociative standard word
where the brackets mean (\ref{op1}).
A {\it super-word} is a word in super-letters.
By 5s every standard word $u$
defines a super-letter $[u].$
\label{sup1}
\end{definition}

Let $D$
be a linearly ordered Abelian additive group.
Suppose that some positive $D$-degrees $d_1,\ldots ,d_n\in D$
are associated to  $x_1,\ldots , x_n.$
We define the degree
of a word to be equal to $m_1d_1+\ldots +m_nd_n$
where $(m_1,\ldots ,m_n)$
is the constitution of the word.
The order and the degree
on the super-letters are defined in the following way:
$[u]>[v]\iff u>v;$\ ${\rm D}([u])={\rm D}(u).$
\addtocounter{nom}{1}
\begin{definition} \rm
A super-letter $[u]$
is called {\it hard in }$H$
provided that its value in $H$
is not a linear combination
of values of super-words of the same degree in less than $[u]$
super-letters
\underline{and $G$-super-words of a lesser  degree}.
\label{tv1}
\end{definition}
\addtocounter{nom}{1}
\begin{definition} \rm
We say that a {\it height} of a super-letter $[u]$
of degree $d$
equals $h=h([u])$
if  $h$
is the smallest number such that: first $p_{uu}$
is a primitive $t$-th root of unity and either $h=t$
or $h=tl^r,$ where $l=$char({\bf k}); and  then the value in $H$
of $[u]^h$
is a linear combination of super-words of degree $hd$
in  less than $[u]$
super-letters \underline{and  $G$-super-words of a lesser degree}.
If there exists no such number  then the height equals infinity.
\label{h1}
\end{definition}

Clearly, if the algebra $H$
is $D$-homogeneous
then one may omit the underlined parts of the above definitions.
\addtocounter{nom}{1}
\begin{theorem}$\!\!([{\bf \ref{pKh1}},\hbox{\rm  Theorem 2}]).$
The set of all values in $H$ of all $G$-super-words $W$
in the hard super-letters $[u_i],$
\begin{equation}
W=g[u_1]^{n_1}[u_2]^{n_2}\cdots [u_m]^{n_m},
\label{BWa}
\end{equation}
where $g\in G,$\ $u_1<u_2<\ldots <u_m,$\ $n_i<h([u_i])$
is a basis of $H.$
\label{BW}
\end{theorem}

In order to find the set of PBW-generators
it is necessary first to include in $U$
the values of all hard super-letters, then for each hard super-letter $[u]$
of a finite height, $h([u])=tl^k,$
to add the values of $[u]^t, [u]^{tl},\ldots [u]^{tl^{(k-1)}},$
and next for each hard super-letter of infinite height such that $p_{uu}$
is a primitive $t$-th
root of unity to add the value of $[u]^t.$

Obviously the set of PBW-generators plays the same role
as the basis of the Lie algebra in the PBW theorem does.
Nevertheless the {\bf k}$[G]$-bimodule
generated by the PBW-generators
is not uniquely defined. It depends on the ordering of
the main generators, the $D$-degree, and under the action of antipode
it transforms to a different bimodule of PBW-generators {\bf k}$[G]S(U).$

Another way to construct PBW-generators is connected with
the M.Kashiwara crystallisation idea
[{\bf \ref{pKas}}, {\bf \ref{pKas1}}].
M.Kashiwara considered the main parameter
of the Drinfeld--Jimbo enveloping algebra  as the temperature
of some physical medium. When the temperature
tends to zero the medium crystallises. By this means the  `crystal'
bases must appear. If we replace $p_{ij}$
with zero then $[u,v]$
turns into $uv,$
while $[u]$
turns into $u.$
\addtocounter{nom}{1}
\begin{lemma}$\!\!${\rm (Bases Crystallisation)}.
Under the above crystallisation the set of PBW-generators
constructed in Theorem $\ref{BW}$
turns into another set of PBW-generators.
\label{kri}
\end{lemma}
{\it Proof}. See [{\bf \ref{pKh1}}, Corollary 1].
\addtocounter{nom}{1}
\begin{lemma}$\!\!${\rm (Super-letters Crystallisation)}.
A super-letter $[u]$
is hard in $H$
if and only if the value of $u$
is not a linear combination of lesser words  of the same degree
\underline{and $G$-words
of a lesser degree}.
\label{kri1}
\end{lemma}
{\it Proof}. See  [{\bf \ref{pKh1}}, Corollary 2].
\addtocounter{nom}{1}
\begin{lemma}
Let $B$
be a set of the super-letters containing $x_1,\ldots ,x_n.$
If each pair $[u], [v]\in B,$ $u>v$
satisfies one of the following conditions

$1)\  [[u][v]]$
is not a standard nonassociative word;

$2)$\  the super-letter $[[u][v]]$
is not hard in $H;$

$3)\  [[u][v]]\in B,$\newline
then the set $B$
includes all hard in $H$
super-letters.
\label{tvp}
\end{lemma}
{\it Proof}. Let $[w]$
be a hard super-letter of minimal degree
such that $[w]\notin B.$
Then $[w]=[[u][v]], u>v$
where $[u],$ $[v]$
are hard super-letters. Indeed, if $[u]$
is not hard then by Lemma \ref{kri1}
we have $u=\sum \alpha _iu_i+S,$
where $u_i<u$
and $D(u_i)=D(u),$ $D(S)<D(u).$
 We have $uv=\sum \alpha _iu_iv+Sv,$
where $u_iv<uv.$
Therefore by Lemma \ref{kri1}, the super-letter  $[w]=[uv]$
can not be hard in $H.$ Contradiction.
Similarly, if $[v]$
is not hard then $v=\sum \alpha _iv_i+S,$ $v_i<v,$ $D(v_i)=D(v),$ $ D(S)<D(v).$
Therefore $uv=\sum \alpha _iuv_i+uS,$ $ uv_i<uv,$
and again $[w]$ can not be hard.

Thus, according to the choice of $[w],$ we get  $[u], [v]\in B.$
Since this pair satisfies neither condition
 1)  nor 2),  the condition 3), $[uv]\in B,$
holds. \hfill $\Box $
\addtocounter{nom}{1}
\begin{lemma}
If ${\bf T}\in H$
is a skew primitive element then
$$
{\bf T}=[u]^h+\sum \alpha _iW_i+\sum \beta _jg_jW_j^{\prime }
$$
where $[u]$
is a hard super-letter, $W_i$
are basis super-words in super-letters less than $[u],$
$D(W_i)=hD([u]),$ $D(W_j^{\prime })<hD([u]).$
Here if $p_{uu}$
is not a root of unity then $h=1;$
if $p_{uu}$
is a primitive $t$-th root of unity then $h=1,$
or $h=t,$
or $h=tl^k,$
where $l$
is the characteristic.
\label{prim}
\end{lemma}
{\it Proof}.
Consider an expansion of $\bf T$
in terms of the basis  (\ref{BWa})
\begin{equation}
{\bf T}=\alpha gU+\sum _{i=1}^k\gamma _i g_iW_i+W^{\prime },
\ \ \alpha \neq 0,
\label{u1c}
\end{equation}
where $gU, g_iW_i$
are different basis elements of
maximal degree, and $U$
is one of the biggest words among $U, W_i$
with respect to the lexicographic ordering
of words in the super-letters. On basis expansion of tensors, the element
$\Delta ({\bf T})-{\bf T}\otimes 1-g_t\otimes {\bf T}$
has only one tensor of the form $gU\otimes \ldots \ $
and this tensor equals $gU\otimes \alpha (g-1).$
Therefore $g=1$
and one may apply [{\bf \ref{pKh1}}, Lemma 13]. \hfill $\Box $

\centerline{5. \large \it Groebner--Shirshov relations systems}

\

\vspace{-1mm}
\stepcounter{par}
\noindent Let  $x_1,\ldots ,x_n$
be variables that have  positive degrees $d_1,\ldots ,d_n\in D.$
Recall that a {\it Hall ordering} of words in $x_1,\ldots ,x_n$
is an order when the words are compared  firstly by the degree
and then words of the same degree are
compared by means of the  lexicographic ordering .
Consider a set of relations
\begin{equation}
w_i=f_i,\ \ \ i\in I,
\label{Gri}
\end{equation}
where $w_i$
is a word and $f_i$
is a linear combination of Hall lesser words.
The system  (\ref{Gri}) is said to be
{\it closed under compositions} or a
{\it Groebner--Shirshov relations system} if first none of $w_i$
contains $w_j,$ $i\neq j\in I$
as a sub-word, and then for each pair of words
$w_k,$
$w_j$
such that some non-empty terminal of $w_k$
coincides with an onset of $w_j,$
that is $w_k=w_k^{\prime }v,$ $w_j=vw_j^{\prime },$
the difference (a {\it composition}) $f_kw_j^{\prime }-w_k^{\prime }f_j$
can be reduced to zero {\it in the free algebra} by means  of
a sequence of one sided substitutions
$w_i\rightarrow f_i,$
$i\in I.$
\addtocounter{nom}{1}
\begin{lemma}$\!\!${\rm
(Diamond Lemma} $[{\bf \ref{pBer}}, {\bf \ref{pBok}},
{\bf \ref{pSh2}}]).$
If the system  $(\ref{Gri})$
is closed under compositions then the
words that have none of $w_i$
as sub-words form  a basis of the algebra $H$
defined by $(\ref{Gri}).$
\label{SBB}
\end{lemma}

If none of the words $w_i$
has sub-words $w_j,$ $j\neq i,$
then the converse statement is valid as well.
Indeed, any composition by means of substitutions $w_i\rightarrow f_i$
can be reduced to a linear combination of words
that have no sub-words $w_i.$
Since $f_iw_j^{\prime }-w_i^{\prime }f_j=$
$(f_i-w_i)w_j^{\prime }-w_i^{\prime }(f_j-w_j),$
this linear combination equals zero in $H.$
Therefore all the coefficients have to be  zero.

Since Bases Crystallisation Lemma provides
the basis that consists of words,
the above  note gives a way to construct the Groebner--Shirshov
relations system for any quantum enveloping algebra.

Let $H$
be a character Hopf algebra generated by  skew primitive semi-invariants
$a_1,\ldots ,a_n$
(or a braided bigraded Hopf algebra generated by graiding homogeneous  primitive elements
$a_1,\ldots ,a_n$)
and let $x_1,\ldots ,x_n$
be the related quantum variables. A non-hard in $H$
super-letter $[w]$
is referred to as a {\it minimal } one if first $w$
has no proper standard sub-words that define non-hard super-letters,
and then $w$
has no sub-words $u^h,$
where $[u]$
is a hard super-letter of the height $h.$

By the Super-letters Crystallisation Lemma, for every minimal
non-hard in $H$
super-letter $[w]$
we may write a relation in $H$
\begin{equation}
w=\sum \alpha _iw_i+\sum \beta _jg_jw_j,
\label{Gr1}
\end{equation}
where $w_j,$ $w_i<w$
in the Hall sense, $D(w_i)=D(w),$ $D(w_j)<D(w).$
In the same way if $[u]$
is a hard in $H$
super-letter of a finite height $h$
then
\begin{equation}
u^h=\sum \alpha _iu_i+\sum \beta _jg_ju_j,
\label{Gr2}
\end{equation}
where $u_j,$ $u_i<u^h$
in the Hall sense, $D(u_i)=hD(u),$ $D(u_j)<hD(u).$
The  relations (\ref{AIc}) and the group operation
provide the relations
\begin{equation}
x_ig=\chi ^{x_i}(g)gx_i,\ \ \ g_1g_2=g_3.
\label{Gr3}
\end{equation}
\addtocounter{nom}{1}
\begin{theorem}
The set of relations $(\ref{Gr1}),$ $(\ref{Gr2}),\,$ and $(\ref{Gr3})$
forms a Groebner--Shirshov  system that defines $H.$
The basis determined by this system in Diamond Lemma
coincides with the crystal basis.
\label{Gr}
\end{theorem}
{\it Proof}. The property 4s implies that none of the left hand sides
of  (\ref{Gr1}), (\ref{Gr2}), (\ref{Gr3})
contains another one as a sub-word.
Therefore by the Bases Crystallisation Lemma it is sufficient
to show that the set of all words $c$
determined in the Diamond Lemma
coincides with the crystal basis. By 3s we have
$c=u_1^{n_1}u_2^{n_2}\cdots  u_k^{n_k},$
where $u_1<\ldots <u_k$
is a sequence of standard words.
Every word  $u_i$
define a hard super-letter $[u_i]$ since in the opposite case $u_i,$
and therefore  $c,$
contains a sub-word $w$
that defines a minimal non-hard super-letter $[w]$.
In the same way $n_i$
does not exceed the height of $[u_i].$  \hfill $\Box $
\addtocounter{nom}{1}
\begin{lemma}
In terms of Lemma $\ref{tvp}$
the set of all super-letters $[[u][v]]$
that satisfy the condition $2)$
contains all minimal non-hard super-letters, but
non-hard generators $x_i.$
\label{min}
\end{lemma}
{\it Proof}. If $[w]$
is a minimal non-hard super-letter then $[w]=[[u][v]],$
where $[u],$ $[v]$
are hard super-letters. By Lemma \ref{tvp} we have  $[u], [v]\in B,$
while $[[u][v]]$
neither satisfies  1) nor 3).  \hfill $\Box $

\centerline {6.\large  \it Quantification with constants}

\

\vspace{-1mm}
\stepcounter{par}
\noindent
By means of the Diamond Lemma in some instances
the investigation of a quantification with constants
can be reduced to one of a simple quantification.

Let $H_1=\langle x_1,\ldots ,x_k||F_1\rangle $
be a character Hopf algebra defined by the quantum variables $x_1,\ldots ,x_k$
and the grading homogeneous relations $\{ f=0:f\in F_1\} ,$
while  $H_2=\langle x_{k+1},\ldots ,x_{n}||F_2\rangle $
is a character Hopf algebra defined by the quantum variables
$x_{k+1},\ldots ,x_n$
and  the grading homogeneous relations $\{ h=0:h\in F_2\} .$
Consider the algebra $H=\langle x_1,\ldots ,x_n||F_1, F_2, F_3\rangle ,$
where $F_3$
is the following system of relations with constants
\begin{equation}
[x_i,x_j]=\alpha _{ij}(1-g_ig_j), \ \ \ 1\leq i\leq k<j\leq n.
\label{jim}
\end{equation}
If the  conditions below are met then the character Hopf
algebra structure on $H$
is uniquely determined:
\begin{equation}
p_{ij}p_{ji}=1, \ 1\leq i\leq k<j\leq n;\ \
\chi ^{x_i}\chi ^{x_j}\neq 1 \Longrightarrow \alpha _{ij}=0.
\label{jim1}
\end{equation}
Indeed, in this case the difference $w_{ij}$
between the left and right hand sides of (\ref{jim})
is a skew primitive semi-invariant of the free enveloping algebra
$G\langle x_1,\ldots ,x_n\rangle .$
Consider the ideals of relations $I_1=$id$(F_1)$
and $I_2=$id$(F_2)$
of $H_1$
and $H_2$
respectively. They are, in the present context, Hopf ideals of
$G\langle x_1, \ldots , x_k\rangle $ and $G\langle x_{k+1}, \ldots , x_n\rangle ,$
respectively. Therefore $V=I_1+I_2+$ $\sum ${\bf k}$w_{ij}$
is an antipode stable coideal of $G\langle X\rangle .$
Consequently the ideal generated by $V$
is a Hopf ideal. It remains to note that
this ideal is generated in $G\langle X\rangle $
by $w_{ij}$
and $F_1, F_2.$
\addtocounter{nom}{1}
\begin{lemma}
Every hard in  $H$
super-letter belongs to either $H_1$ or $H_2,$
and it is hard in the related algebra.
\label{ha}
\end{lemma}
{\it Proof}. If a standard word contains at least one
of the letters $x_i,$ $i\leq k$
then it has to start with one of them (see s2).
If this word contains a letter $x_j,$ $j>k$
then it has a sub word of the form $x_ix_j,$ $i\leq k<j.$
Therefore by Lemma \ref{kri1} and relations
(\ref{jim}) this word defines a non-hard super-letter.\hfill $\Box $

The converse statement is not universally true.
In order to formulate the necessary and sufficient conditions
let us define partial skew derivatives:
\begin{eqnarray}
& & (x_j)_i^{\prime }=
 (x_i)_j^{\prime }=\alpha _{ij}(1-g_ig_j),\ \  \ \
i\leq k<j;\hfill \nonumber \\
& &(v\cdot w)_i^{\prime }=
 (v)_i^{\prime }\cdot w+p(x_i,v)v\cdot (w)_i^{\prime },\
i\leq k,\ v,w\in G\langle x_{k+1},
\ldots ,x_n\rangle ;\hfill \nonumber \\
& & (u\cdot v)_j^{\prime }=p(v,x_j) (u)_j^{\prime }\cdot v+
u\cdot  (v)_j^{\prime },\
j>k,\ u,v\in G\langle x_1,\ldots ,x_k\rangle .
\label{paar}
\end{eqnarray}
\addtocounter{nom}{1}
\begin{lemma}
All hard in $H_1$
or $H_2$
super-letters are hard in $H$
if and only if $ (h)_i^{\prime }=0$
in $H_2$
for all $i\leq k,$ $h\in F_2$
and $ (f)_j^{\prime }=0$
in $H_1$
for all $j>k,$ $f\in F_1.$
If these conditions are met  then
 \begin{equation}
              H\cong H_2\otimes _{{\bf k}[G]}H_1
 \label{new}
 \end{equation}
as {\bf k}$[G]$-bimodules and
the space generated by the
skew primitive elements of $H$
equals the sum of these spaces for $H_1$
and $H_2.$
\label{raz}
\end{lemma}

{\it Proof}.  By (\ref{dif1}) and  (\ref{paar})
the following equalities are valid in $H:$
\begin{equation}
0=[x_i,h]= (h)_i^{\prime };\ \
0=[f,x_j]= (f)_j^{\prime },\ i\leq k<j.
\label{prd}
\end{equation}
If all hard in $H_1$
or $H_2$
super-letters are hard in $H$
then $H_1,$ $H_2$
are sub-algebras of $H.$
So (\ref{prd}) proves the necessity of the lemma conditions.

Conversely. Let us consider an algebra $R$
defined by the generators $g\in G,$ $x_1,\ldots , x_n$
and  the relations (\ref{Gr3}), (\ref{jim}).
Evidently this system is closed under the compositions.
Therefore by Diamond Lemma
the set of words $gvw$
forms a basis of $R$
where  $g\in G;$  $v$
is a word in $x_j,$ $j>k;$
and $w$
is a word in $x_i,$ $i\leq k.$
In other words $R$
as a bimodule over  {\bf k}$[G]$
has a decomposition
\begin{equation}
R=G\langle x_{k+1},\ldots ,x_n\rangle \otimes _{\hbox{\bf k}[G]}
G\langle x_1,\ldots ,x_k\rangle .
\label{prd1}
\end{equation}

Let us show that the two sided ideal of $R$
generated by $F_2$
coincides with the right ideal $I_2R=$
$I_2\otimes _{\hbox{\bf k}[G]} G\langle x_1,\ldots ,x_k\rangle .$
It will suffice to show that $I_2R$
admits left multiplication by $x_i,$ $i\leq k.$
If  $v$
is a word in $x_{k+1}, \ldots ,x_k,$ $h\in F_2,$ $r\in R$
then $x_ivhr=[x_i, vh]r+p(x_i,vh)vhx_ir.$
The second term belongs to $I_2R,$
while the first one can be rewritten by  (\ref{dif1}):
$[x_i,v]h+p(x_i,v)v[x_i,h].$
Both of these addends belong to $I_2R$
since  $[x_i,v]=$ $(v)_i^{\prime }\in $
$G\langle x_{k+1},\ldots ,x_n\rangle $
and $[x_i,h]=$ $(h)_i^{\prime }\in I_2.$

Furthermore, consider a quotient algebra $R_1=R/I_2R:$
$$
R_1=(G\langle x_{k+1},\ldots ,x_n\rangle \otimes _{\hbox{\bf k}[G]}
G\langle x_1,\ldots ,x_k\rangle )/
(I_2\otimes _{\hbox{\bf k}[G]}G\langle x_1,\ldots ,x_k\rangle )=
$$
$$
H_2\otimes _{\hbox{\bf k}[G]}G\langle x_1,\ldots ,x_k\rangle ,
$$
where the equality means the natural isomorphism of
{\bf k}$[G]$-bimodules.

Along similar lines,  the left ideal $R_1I_1=$
$H_2\otimes _{\hbox{\bf k}[G]}I_1 $
of this quotient algebra coincides with the two sided ideal
generated by $F_1.$
Therefore
$$
H=R_1/R_1I_1=
H_2\otimes _{\hbox{\bf k}[G]}G\langle x_1,\ldots ,x_k\rangle /
H_2\otimes _{\hbox{\bf k}[G]}I_1=H_2\otimes _{\hbox{\bf k}[G]}H_1.
$$
Thus the monotonous restricted $G$-words
in hard in $H_1$
or $H_2$
super-letters form a basis of $H.$
This, in particular,  proves the first statement.

Now let $T=\sum \alpha _tg_tV_tW_t$
be the basis decomposition of a skew primitive element,
$g_t\in G,$ $V_t\in H_2,$ $W_t\in H_1,$ $\alpha _t\neq 0.$
We have to show that for each $t$
one of the super-words $V_t$
or $W_t$
is empty. Suppose that it is not so. Among the
addends with non-empty $V_t,$ $W_t$
we choose the largest one in the Hall sense, say $g_sV_sW_s.$
Under the basis decomposition of $\Delta (T)-T\otimes 1-g(T)\otimes T$
the term $\alpha _sg_sg(V_s)W_s\otimes g_sV_s$
appears and cannot be cancelled  with other.
Indeed, since the coproduct is
homogeneous (see [{\bf \ref{pKh1}}, Lemma 9])
and since under the basis decomposition the super-words are decreased
(see [{\bf \ref{pKh1}}, Lemma 7]) the product
$\alpha _s(g_s\otimes g_s)\Delta (V_s)\Delta (W_s)$
has the only term of the above type.
By the same reasons $\alpha _t(g_t\otimes g_t)\Delta (V_t)\Delta (W_t)$
has a term of the above type only if $V_t\geq V_s$
and $W_t\geq W_s$
with respect to the  Hall ordering  of the set of all super-words.
However, by the choise of $s,$
we have $D(V_sW_s)\geq D(V_tW_t).$
Hence $D(V_t)=D(V_s)$
and $D(W_t)=D(W_s).$
In particular  $V_t$
is not a proper onset of  $V_s.$
Therefore $V_t=V_s$
since otherwise the inequality $V_t>V_s$
yields a contradiction $V_tW_t>V_sW_s$.
The inequality $W_t>W_s$
get the same contradiction. Therefore  $V_t=V_s$
and $W_t=W_s,$
in which case $g_tg(V_t)W_t\otimes g_tV_t=$ $g_sg(V_s)W_s\otimes g_sV_s.$
Thus $g_t=g_s$
and $t=s.$ \hfill $\Box $

\

\centerline{7.\large  \it Quantification of the classical series}

\

\vspace{-1mm}
\stepcounter{par}
\noindent
In this section we apply the above general results to the
infinite series $A_n,$ $B_n,$ $C_n,$ $D_n$
of nilpotent Lie algebras
defined by the Serre relations (\ref{rel}). Let $\frak g$ be any such Lie algebra.
\addtocounter{nom}{1}
\begin{lemma}
If a standard word $u$
has no sub words of the type
\begin{equation}
x_i^sx_jx_i^m,\ \hbox{where}\  s+m=1-a_{ij}
\label{vor}
\end{equation}
then $[u]$
is a hard in $U_P({\frak g})$
super-letter.
\label{tver}
\end{lemma}
{\it Proof}. Let $R$
be defined by the generators $x_1,\ldots ,x_n$
and the relations
\begin{equation}
x_i^sx_jx_i^m=0,\ \hbox{ where }\ s+m=1-a_{ij}.
\label{R2}
\end{equation}
Clearly (\ref{R2}) implies (\ref{rel}) with
the skew commutator  in place of
the Lie operation. Therefore  $R$
is a homomorphic image of $U_P({\frak g}).$
The system  (\ref{R2}) is closed under compositions
since a composition of monomial relations always has the form $0=0.$

Let  $u$
have no sub-words (\ref{vor}). If $[u]$
is  not hard then, by the Super-letters Crystallisation Lemma, $u$
is a linear combination of lesser words in $U_P({\frak g}).$
Therefore $u$
is a linear combination of lesser words in $R$
as well. This contradicts  the fact that $u$
belongs to the Groebner--Shirshov  basis of $R,$
since every word either belongs to this basis or equals zero in $R$.  \hfill $\Box $

\

{\bf Theorem $\bf A_n$}. {\it
Suppose that ${\frak g}$
is of the type $A_n,$
and $p_{ii}\neq -1.$
Denote by $B$
the set of the super-letters given below:
\begin{equation}
[u_{km}]\stackrel{df}{=}[x_kx_{k+1}\ldots x_m],
\ \ \ 1\leq k\leq m\leq n.
\label{ANb}
\end{equation}
The following statements are valid.

$1.$ The values of $[u_{km}]$
in $U_P({\frak g})$
form a PBW-generators set.

$2.$ Each of the super-letters $(\ref{ANb})$ has infinite height in
$U_P({\frak g}).$

$3.$ The values of all non-hard in $U_P({\frak g})$
super-letters equal zero.

$4.$ The following relations  with $(\ref{Gr3})$ form
the Groebner--Shirshov relations system that determines
the crystal basis of $U_P({\frak g}):$
\begin{equation}
\begin{matrix}
[u_0]\stackrel{df}{=}[x_kx_m]=0,\hfill &1\leq k<m-1<n;\hfill \cr
[u_1]\stackrel{df}{=}[x_kx_{k+1}\ldots x_mx_{k+1}]=0,
\hfill &1\leq k<m\leq n;\hfill \cr
[u_2]\stackrel{df}{=}[x_kx_{k+1}\ldots
x_mx_kx_{k+1}\ldots x_{m+1}]=0,\hfill &
1\leq k\leq m<n.\hfill \end{matrix}
\label{GAn}
\end{equation}

$5.$ If $p_{11}\neq 1$
then the generators $x_i,$
the constants $1-g,$ $g\in G,$
and, in the case that $p_{11}$
is a primitive $t$-th root of 1, the elements $x_i^t, x_i^{tl^k}$
form a basis of ${\frak g}_P=L(U_P({\frak g})).$
Here $l$
is the characteristic of the ground field.

$6.$ If  $p_{11}=1$
then the elements $(\ref{ANb})$
and, in the case $l>0,$
their $l^k$-th
powers, together with $1-g,$ $g\in G$
form a basis of ${\frak g}_P.$}

\

By Corollary \ref{sled} the relations (\ref{rel})
with a Cartan matrix $A$
of type $A_n$
admit a quantification if and only if
\begin{equation}
p_{ii}=p_{11},\ \
p_{ii+1}p_{i+1i}=p_{11}^{-1};\ \
p_{ij}p_{ji}=1,\ i-j>1.
\label{suxa}
\end{equation}
In this case the quantified relations (\ref{rel}) take up the form
\begin{eqnarray}
& & x_ix_{i+1}^2=p_{i i+1}(1+p_{i+1 i+1})x_{i+1}x_ix_{i+1}-
p_{i i+1}^2p_{i+1 i+1}x_{i+1}^2x_i,
\label{rA1} \\
& & x_i^2x_{i+1}=p_{i i+1}(1+p_{ii})x_ix_{i+1}x_i-
p_{i i+1}^2p_{ii}x_{i+1}x_i^2,
\label{rA2} \\
& & x_ix_j=p_{ij}x_jx_i,\ \ \ \ i-j>1.
\label{rA3}
\end{eqnarray}

Let us introduce a congruence   $u\equiv _kv$
on $G\langle X\rangle .$
This congruence means that the value of  $u-v$
in $U_P^b({\frak g})$
belongs to the subspace generated by values of all
words with the initial letters  $x_i, i\geq k.$
Clearly, this congruence admits right multiplication
by arbitrary polynomials as well as left
multiplication by the independent of $x_{k-1}$
ones (see (\ref{rA3})). For example,
by (\ref{rA1}) and (\ref{rA2}) we have
\begin{equation}
x_ix_{i+1}^2\equiv _{i+1}0;\ \ \
x_ix_{i+1}x_i\equiv _{i+1}\alpha x_i^2x_{i+1},\ \ \ \alpha \neq 0.
\label{cong1a}
\end{equation}
\addtocounter{nom}{1}
\begin{lemma}
If $y=x_i,\ m+1\neq i>k$
or $y=x_i^2,\ m+1=i>k$
then
\begin{equation}
u_{km}y\equiv _{k+1}0.
\label{cong4a}
\end{equation}
\label{ong4a}
\end{lemma}
{\it Proof}. Let $y=x_{m+1}^2,$ $m+1>k.$
By (\ref{cong1a}) and (\ref{rA3}) we have that
$u_{km}y=$ $u_{k\, m-1}\underline{x_mx_{m+1}^2}\equiv _{m+1}0.$
If $y=x_i$ and  $m+1\neq i>k$
then we get $u_{km}y=$
$\alpha u_{k\, i-1}\underline{x_ix_{i+1}x_i}$
$u_{i+2m}\equiv _{i+1}
\beta \underline{u_{k\, i-1}x_i^2}u_{i+1m}\equiv _{k+1}0$
by the above case. \hfill $\Box $
\addtocounter{nom}{1}
\begin{lemma}
The brackets in $[u_{km}]$
are left-ordered, $[u_{km}]=[x_k[u_{k+1m}]].$
\label{an1}
\end{lemma}
{\it Proof}. The  statement immediately follows from  the properties  6s and 2s. \hfill $\Box $
\addtocounter{nom}{1}
\begin{lemma} If a nonassociative word $[[u_{km}][u_{rs}]]$
is standard then
$ k=m\leq r;$ or  $r=k+1,$  $m\geq s;$ or  $ r=k,$ $m<s.$
\label{an2}
\end{lemma}
{\it Proof}. By definition, $u_{km}>u_{rs}$
if and only if either $k<r;$ or $k=r,$ $m<s.$
If $k=m$
then $u_{km}=x_k$
and $m\leq r.$
If $k\neq m$
then $[u_{km}]=[x_k[u_{k+1m}]].$
Therefore $u_{k+1m}\leq u_{rs},$
i.e. either $k+1>r;$
or $k+1=r$
and $m\geq s.$
The former case contradicts $k<r$
while the latter one does $k=r.$
Thus only the possibilities set in the lemma remain. \hfill $\Box $
\addtocounter{nom}{1}
\begin{lemma} If $[w]=[[u_{km}][u_{rs}]],$ $n\geq 1$
is a standard nonassociative word then the constitution of $[w]^h$
does not equal the constitution of any super-word in less than $[w]$
super-letters from $B.$
\label{man}
\end{lemma}
{\it Proof}. The inequalities at the last column of the  following
tableaux are valid  for all $[u]\in B$
that are less than the super-letters located on the same row,
where as above  deg$_i(u)$
means  the degree of $u$
in $x_i.$
\begin{equation}
\begin{matrix}
[x_ku_{k+1s}] \hfill & &\hbox{deg}_k(u)\leq \hbox{deg}_{s+1}(u); \hfill \cr
[x_ku_{rs}], \hfill & k\leq r\neq k+1 \hfill &
\hbox{deg}_k(u)\leq \hbox{deg}_{k+1}(u); \hfill \cr
[u_{km}u_{k+1s}], \hfill & m\geq s \hfill &
\hbox{deg}_k(u)\leq \hbox{deg}_{m+1}(u); \hfill \cr
[u_{km}u_{ks}], \hfill & m<s \hfill &
\hbox{deg}_k(u)\leq \hbox{deg}_{m+1}(u). \hfill \end{matrix}
\label{maf}
\end{equation}
If all super-letters of a super-word $U$
satisfy one of these  inequalities then $U$
does as well. Clearly, no one of the super-letters in the first column
satisfies the degree inequality on the same row. Finally, by Lemma  \ref{an2} the first column
contains all standard nonassociative words
of the type $[[u_{km}][u_{rs}]].$ \hfill $\Box $
\addtocounter{nom}{1}
\begin{lemma}
If $p_{11}\neq 1$
then the values of $[u_{km}]^h,$ $k<m,$ $h\geq 1$
are not skew primitive, in particular they are non-zero.
\label{alpr}
\end{lemma}
{\it Proof}. The sub-algebra generated by $x_2,\ldots x_n$
is defined by the Cartan matrix of the type $A_{n-1}.$
This allows us to use induction on $n.$
If $n=1$
then the lemma is correct in the sense that $[u_{km}]^h=x_1^h\neq 0.$

Let $n>1.$
If $k>1$
then we may use the inductive supposition directly.
Consider the decomposition $\Delta ([u_{1m}])=\sum u^{(1)}\otimes u^{(2)}.$
Since
\begin{equation}
[u_{1m}]=x_1[u_{2 m}]-p(x_1,u_{2 m})[u_{2 m}]x_1,
\label{pov1}
\end{equation}
we have
$$
\Delta ([u_{1m}])=(x_1\otimes 1+g_1\otimes x_1)
\Delta ([u_{2 m}])-
$$
\begin{equation}
p(x_1,u_{2 m})\Delta ([u_{2 m}])
(x_1\otimes 1+g_1\otimes x_1).
\label{pov2}
\end{equation}
Therefore the sum of all tensors $u^{(1)}\otimes u^{(2)}$
with deg$_1(u^{(2)})=1,$ deg$_k(u^{(2)})=0,$ $k>1$
has the form $\varepsilon g_1[u_{2 m}]\otimes x_1,$
where  $\varepsilon =$
$1-p(x_1,u_{2 m})p(u_{2 m},x_1)$
since $[u_{2 m}]g_1=p(u_{2 m},x_1)g_1[u_{2 m}].$
By (\ref{suxa}) we have $p_{ij}p_{ji}=1$
for $i-1>j.$
Therefore $\varepsilon =1-p_{1 2}p_{2 1}=$ $1-p_{11}^{-1}\neq 0.$

This implies that in the decomposition
$\Delta ([u_{1m}]^h)=\sum v^{(1)}\otimes v^{(2)}$
the sum of all tensors $v^{(1)}\otimes v^{(2)}$
with deg$_1(v^{(2)})=h,$ deg$_k(v^{(2)})=0,$ $k>1$
equals $\varepsilon ^h[u_{2 m}]^h\otimes x_1^h.$
Thus $[u_{1m}]^h$
is not skew primitive in $U_P({\frak g}).$ \hfill $\Box $

{\it Proof} of Theorem ${\bf A_n}.$
Let us show firstly that $B$
satisfies  the conditions of Lemma \ref{tvp}.
By the Super-letter Crystallisation Lemma $[w]=[[u_{km}][u_{rs}]]$
is non-hard if the value of $u_{km}u_{rs}$
is a linear combination of lesser words. For $k=m,$ $r=k+1$
we have $[w]=[u_{ks}]\in B.$
If $k=m,$ $r>k+1$
then the word $x_ku_{rs}$
can be diminished by (\ref{rA2}) or (\ref{rA3}).
If $k\neq m$
then by Lemma \ref{an2} the word $u_{km}u_{rs}$
has a sub-word of the type $u_1$
or $u_2.$
Thus we need show only that the values in $U_P({\frak g})$
of $u_1$ and $u_2$ are linear combinations of lesser words.

The word $u_1$
has such a representation by Lemma \ref{ong4a}. Consider the word $u_2.$
Let us show by downward induction on $k$
that
\begin{equation}
u_{k m}u_{k\, m+1} \equiv _{k+1}\gamma u_{k\, m+1}u_{km},\ \ \ \gamma \neq 0.
\label{ume}
\end{equation}
If $k=m$
then one may use (\ref{rA2}) with $i=k.$
Let $k<m.$
Let us transpose the second letter $x_k$
of $u_2$
as far to the left  as possible by (\ref{rA3}). We get
$$
u_2=\alpha \underline{x_kx_{k+1}x_k}x_{k+2}\cdots
x_mx_{k+1}\cdots x_{m+1}, \ \ \ \alpha \neq 0.
$$
By (\ref{rA2}) we have
$$
u_2\equiv _{k+1}\beta x_k^2(x_{k+1}x_{k+2}\cdots x_mx_{k+1}\cdots
x_{m+1}),\ \ \beta \neq 0.
$$
Let us apply the inductive supposition to the word in the parentheses.
Since  $x_i,$ $i>k+1$
commutes with $x_k^2$
according to the formulae (\ref{rA3}), we get
$$
u_2\equiv _{k+1}\gamma \underline{x_k^2x_{k+1}}x_{k+2}\cdots
x_{m+1}x_{k+1}\cdots x_m.
$$
Now it remains to replace the underlined sub-word according to
(\ref{rA2}) and then to transpose the second letter $x_k$
to its former position by (\ref{rA3}).

Note that for the diminishing of $u_1,$ $u_2$
we did not use, and we could not use, the relation $[x_{n-1}x_n^2]=0$
since  deg$_n(u_1)\leq 1,$ deg$_n(u_2)\leq 1.$

Thus $B$ satisfies the conditions of Lemma \ref{tvp}.
Since none of  $[u_{km}]$ has sub-words (\ref{vor}),
Lemmas  \ref{tver} and \ref{tvp} show
that the first statement is correct.

If  $[u_{km}]$
has a finite height $h$
then the value of the polynomial $[u_{km}]^h$
in $U_P({\frak g})$
is a linear combination of words in hard super-letters
that are less than $[u_{km}].$
However by Lemma \ref{man}  this linear combination is trivial,
$[u_{km}]^h=0,$
since the defining relations are homogeneous. By Lemma \ref{alpr}
the second statement is correct for $p_{11}\neq 1$.

Similarly consider the skew primitive elements. Since both
the defining relations and the coproduct are homogeneous,
all the homogeneous components of a skew primitive element
are skew primitive itself. Therefore it remains to describe
all skew primitive elements homogeneous  in each $x_i.$ Let $T$
be such an element. By Lemma \ref{prim} we have
$$
T=[u]^h+\sum \alpha _iW_i,
$$
where $[u]$
is a hard super-letter, $u=u_{km},$
and $W_i$
are  super-words in less than $[u]$
super-letters from $B$. By the homogeneity all $W_i$
have the same constitution as $[u_{km}]^h$
does. However by Lemma  \ref{man} there exist no such  super-words.
This means that the only case possible is  $T=[u_{km}]^h.$
Thus, by Lemma \ref{alpr} the fifth statement is valid as well.

If $p_{11}=1$
then $p_{ij}p_{ji}=p_{ii}=1$
for all $i,j.$
So we are under the conditions of Example 1,
that is $U_P^b({\frak g})$
is the universal enveloping algebra of the colour Lie algebra
${\frak g}^{col}.$
Further, $[u_{km}]\in {\frak g}^{col}$
and $[u_{km}]$
are linearly independent in ${\frak g}^{col}$
since they are hard super-letters and no one of them
can be a linear combination of the lesser ones.
Let us complete $B$
to a homogeneous basis $B^{\prime }$
of ${\frak g}^{col}.$
Then by the PBW theorem for the colour Lie algebras the
products $b_1^{n_1}\cdots b_k^{n_k},$ $b_1<\ldots <b_k$
form a basis of
$U({\frak g}^{col})=$ $U_P^b({\frak g}).$
However, the monotonous restricted words in $B$
form a basis of $U_P^b({\frak g})$
also. Thus $B^{\prime }=B$
and all hard super-letters have the infinite height.

In particular, we get that the second statement is valid in complete
extent. Moreover, if $p_{11}=1$
then $p(u_{km},u_{km})=1,$
thus for $l=0$ all homogeneous skew primitive
elements became exhausted by $[u_{km}],$
while for $l>0$
the powers $[u_{km}]^{l^k}$
are added to them (of course, here $l\neq 2$
since $-1\neq p_{ii}=1$).

So we have  proved all statements, but the third and fourth ones.
These statements will follow Theorem \ref{Gr}
and Lemma \ref{min} if we prove that all non-hard
super-letters  $[[u_{km}][u_{rs}]]$
equal zero in $U_P({\frak g}).$
By the homogeneous definition, $[[u_{km}][u_{rs}]]$
is a linear
combination of super-words in lesser hard super-letters.
However,  by Lemma \ref{man}, there exist no such super-words
of the same constitution. Therefore, by the homogeneity,
the above linear combination equals zero. \hfill $\Box $

{\bf Theorem $\bf B_n$}. {\it
Let ${\frak g}$
be of the type $B_n,$
and $p_{ii}\neq -1,$ $1\leq i<n,$ $p_{nn}^{[3]}\neq 0.$
Denote by $B$
the set of the super-letters given below:
\begin{equation}
\begin{matrix}
[u_{km}]\stackrel{df}{=}[x_kx_{k+1}\ldots x_m],\hfill &
1\leq k\leq m\leq n;\hfill \cr
[w_{km}]\stackrel{df}{=}[x_kx_{k+1}\ldots x_n\cdot x_n\ldots x_m],\hfill &
1\leq k<m\leq n.\hfill \end{matrix}
\label{BN1}
\end{equation}
The following statements are valid.

$1.$ The values of $(\ref{BN1})$
in $U_P({\frak g})$
form the PBW-generators set.

$2.$ Every super-letter   $[u]\in B$
has infinite height in $U_P({\frak g}).$

$3.$ The relations $(\ref{Gr3})$ with the following ones
form a Groebner--Shirshov
system that determines the crystal basis of $U_P({\frak g}).$
\begin{equation}
\begin{matrix}
[u_0]\stackrel{df}{=}[x_kx_m]=0,\hfill &
1\leq k<m-1<n;\hfill \cr
[u_1]\stackrel{df}{=}
[u_{km}x_{k+1}]=0,\hfill &
1\leq k<m\leq n,\  k\neq n-1;\hfill \cr
[u_2]\stackrel{df}{=}
[u_{km}u_{k\, m+1}]=0,\hfill &
1\leq k\leq m<n;\hfill \cr
[u_3]\stackrel{df}{=}
[w_{km}x_{k+1}]=0,\hfill &
1\leq k<m\leq n,\  k\neq m-2;\hfill \cr
[u_4]\stackrel{df}{=}
[w_{k k+1}x_{k+2}]=0,\hfill &
1\leq k<n-1;\hfill \cr
[u_5]\stackrel{df}{=}
[w_{km}w_{k\, m-1}]=0,\hfill &
1\leq k<m-1\leq n-1;\hfill \cr
[u_6]\stackrel{df}{=}
[u_{kn}^2x_n]=0,\hfill &
1\leq k<n.\hfill \end{matrix}
\label{bGB}
\end{equation}

$4.$
If $p_{11}\neq 1$
then the generators $x_i$
and their powers $x_i^t, x_i^{tl^k},$
such that $p_{ii}$
is a primitive $t$-th root of 1,
together with the constants $1-g,$ $g\in G$
form a basis of ${\frak g}_P=L(U_P({\frak g})).$
Here $l$
is the characteristic of the ground field.

$5.$
If  $p_{nn}=p_{11}=1$
then the elements $(\ref{BN1})$
and,  for $l>0,$
their $l^k$-th
powers, together with $1-g,$ $g\in G$
form a basis of ${\frak g}_P.$
If $p_{nn}=-p_{11}=-1$
then $[u_{kn}]^2,$ $[u_{kn}]^{2l^k}$
are added to them}.

\

Recall that in the case $B_n$
the algebra $U_P^b({\frak g})$
is defined by (\ref{rA1}), (\ref{rA2}), (\ref{rA3})
where in (\ref{rA1})  the last relation, $i=n-1,$
is replaced with
\begin{equation}
x_{n-1}x_n^3=p_{n-1 n}p_{nn}^{[3]}x_nx_{n-1}x_n^2-
p_{n-1 n}^2p_{nn}p_{nn}^{[3]}x_n^2x_{n-1}x_n+
p_{n-1 n}^3p_{nn}^3x_n^3x_{n-1}.
\label{rB1}
\end{equation}
By Corollary  \ref{sled} we get the existence conditions
\begin{equation}
p_{ii}=p_{11},\ p_{ii+1}p_{i+1i}=p_{11}^{-1}=p_{nn}^{-2},\
1\leq i\leq n-1;\ p_{ij}p_{ji}=1,\ i-j>1.
\label{sux}
\end{equation}
The relations (\ref{rA1}) and (\ref{rB1}) show that
\begin{equation}
x_ix_{i+1}^2\equiv _{i+1}0,\ \ i<n-1;\ \ \ \ \
x_{n-1}x_n^3\equiv _n0,
\label{cong1}
\end{equation}
while the relations  (\ref{rA2}) imply
\begin{equation}
x_ix_{i+1}x_i\equiv _{i+1}\alpha x_i^2x_{i+1},\ \ \ \alpha \neq 0.
\label{cong2}
\end{equation}
By means of these relations and (\ref{rA3}), (\ref{rB1}) we have
\begin{equation}
x_{n-2}x_{n-1}\underline{x_n^2x_{n-1}x_n}\equiv _{n-1}0.
\label{ng3}
\end{equation}
\addtocounter{nom}{1}
\begin{lemma}
The brackets in $[w_{km}]$
are set by the  recurrence formulae:
\begin{equation}
\begin{matrix}
[w_{km}]=[x_k[w_{k+1 m}]],\hfill &\hbox{if}\ 1\leq k<m-1<n;\hfill \cr
[w_{k k+1}]=[[w_{k k+2}]x_{k+1}],\hfill &\hbox{if}\ 1\leq k<n.\end{matrix}
\label{sb}
\end{equation}
Here by the definition $w_{k\, n+1}=u_{k n}.$
\label{bn1q}
\end{lemma}
{\it Proof}. It is enough to use the property 6s and then 1s and 2s. \hfill $\Box $
\addtocounter{nom}{1}
\begin{lemma} The nonassociative word $[[w_{km}][w_{rs}]]$
is standard only in the following two cases: $1)\ s\geq m>k+1=r;$
$2)\ s<m,$ $r=k.$
\label{bn2}
\end{lemma}
{\it Proof}. If $[[w_{km}][w_{rs}]]$
is standard then $w_{km}>w_{rs}$
and by (\ref{sb}) either $w_{k+1}\leq w_{rs},$
or $m=k+1$
and $x_{k+1}\leq w_{rs}.$
The inequality $w_{km}>w_{rs}$
is correct only in two cases:  $k<r$
or $k=r, m>s.$
We get four possibilities:
$1)\  k<r,$ $\ k<m-1,$ $\ w_{k+1 m}\leq w_{rs};$
$2)\  k<r,$ $\ m=k+1,$ $\ x_{k+1}\leq w_{rs};$
$3)\  k=r,$ $\ m>s,$ $\ k<m-1,$ $\ w_{k+1 m}\leq w_{rs};$
$4)\  k=r,$ $\ m>s,$ $\ m=k+1,$ $\ x_{k+1}\leq w_{rs}.$
Only the first and third ones are consistent since in the second case
$x_{k+1}\leq w_{rs}$
implies $k+1>r,$
while in the fourth case $r<s$
and $k=r<s<m=k+1.$
If now  we decode $w_{k+1 m}\leq w_{rs}$
in the first and third cases, we get  the
two possibilities mentioned in the lemma. \hfill $\Box $
\addtocounter{nom}{1}
\begin{lemma} The nonassociative word $[[u_{km}][w_{rs}]]$
is standard only in the following two cases: $1)\ k=r;$ $2)\ k=m<r.$
\label{bn3}
\end{lemma}
{\it Proof}. The inequality $u_{km}>w_{rs}$
means $k\leq r.$
Since $[u_{km}]=[x_k[u_{k+1 m}]],$
for $k\neq m$
we get $u_{k+1 m}\leq w_{rs},$
so $k+1>r$
and $k=r.$
If $k=m\neq r$
then $x_m>w_{rs}$
and $m<r.$ \hfill $\Box $
\addtocounter{nom}{1}
\begin{lemma} The  nonassociative word $[[w_{km}][u_{rs}]]$
is standard only in the following two cases: $1)\ r=k+1<m;$
$2)\ r=k+1=m=s.$
\label{bn4}
\end{lemma}
{\it Proof}. The inequality $w_{km}>u_{rs}$
implies $r>k.$
If $k<m-1$
then by the first formula (\ref{sb}) we have $w_{k+1 m}\leq u_{rs}$
that is equivalent to $k+1\geq r.$
Therefore  $r=k+1<m.$
If $k=m-1$
then by the second formula (\ref{sb}) we get $x_{k+1}\leq u_{rs},$
i.e. either $k+1>r$
or $k+1=r=s.$
The former case contradicts $r>k$
while the latter one is mentioned in the lemma. \hfill $\Box $
\addtocounter{nom}{1}
\begin{lemma}
If $[u], [v]\in B$
then one of the statements below is correct.

$1)\ [[u][v]]$
is not a standard nonassociative word;

$2)\  uv$
contains a sub-word of one of the types
$u_0, u_1, u_2, u_3, u_4, u_5, u_6;$

$3)\  [[u][v]]\in B.$
\label{bn6}
\end{lemma}
{\it Proof}. The proof results from  Lemmas
\ref{an2}, \ref{bn2}, \ref{bn3}, \ref{bn4}. \hfill $\Box $
\addtocounter{nom}{1}
\begin{lemma} If a super-word $W$
equals one of the super-letters $[u_1]$--$[u_6]$
or $[u_{km}]^h,$ $[w_{km}]^h,$ $h\geq 1$
then its constitution does not equal the constitution
of any super-word in less than $W$
super-letters from $B.$
\label{mbn}
\end{lemma}
{\it Proof}. The proof is akin to Lemma \ref{man}
with the following tableaux:
\begin{equation}
\begin{matrix}
[u_{km}],\ \, [u_{km}x_{k+1}], \ \, [u_{km}u_{k\, m+1}] \hfill &
\hbox{deg}_k(u)\leq \hbox{deg}_{m+1}(u); \hfill \cr
[w_{km}],\ [w_{km}x_{k+1}], \ [w_{km}w_{k\, m-1}] \hfill &
2\hbox{deg}_k(u)\leq \hbox{deg}_{m-1}(u); \hfill \cr
[w_{kk+1}x_{k+2}] \hfill &
\hbox{deg}_k(u)=0; \hfill \cr
[u_{kn}^2x_n]\hfill &
\hbox{deg}_k(u)\leq \hbox{deg}_n(u).\hfill \end{matrix}
\label{mbf}
\end{equation}

\hfill $\Box $

\addtocounter{nom}{1}
\begin{lemma}
If $y=x_i,\ m-1\neq i>k$
or $y=x_i^2,\ m-1=i>k$
then
\begin{equation}
w_{km}y\equiv _{k+1}0.
\label{cong4}
\end{equation}
\label{ong4}
\end{lemma}
{\it Proof}. If $i<m-1$
then by means of (\ref{rA3}) it is possible to permute $y$
to the left beyond $x_n^2$
and use Lemma  \ref{ong4a} with $m^{\prime }=n-1.$
If $y=x_i^2,\ $ $m-1=i>k$
then by the above case, $i<m-1,$ we get
\begin{equation}
w_{km}y=w_{k\, m+1}\underline{x_mx_{m-1}^2}=
\underline{w_{k\, m+1}x_{m-1}}
(\alpha x_mx_{m-1}+\beta x_{m-1}x_m)\equiv _{k+1}0,
\label{sss}
\end{equation}
where for $m=n$
by definition $w_{k\, n+1}=$ $u_{kn},$
and $u_{kn}x_{n-1}\equiv _{n-1}0.$

If $y=x_i, i=m>k$
then for $m=n$
one may use the second equality (\ref{cong1}). For $m<n$
we have $w_{km}y=w_{k\, m+1}y_1$
where $y_1=x_m^2.$
Therefore for  $k<n-1$
we may use (\ref{sss}) with $m+1$
in place of $m.$
For  $k=n-1$
we have  $w_{km}x_n=x_{n-1}x_n^3\equiv _n0.$

Finally, if $y=x_i,\ i>m>k$
then by (\ref{rA3}) we have
$w_{km}y=\alpha w_{k i+1}\underline{x_i x_{i-1}x_i}\cdot v.$
For  $i=n$
one may use (\ref{ng3}), while for  $i<n,$
changing the underlined word according to (\ref{rA1}),
we may use the above considered cases: $m^{\prime }-1=i^{\prime },$
where $m^{\prime }=i+1,$ $i^{\prime }=i;$
and $i^{\prime }<m^{\prime }-1,$
where $m^{\prime }=i+1,$ $i^{\prime }=i-1.$ \hfill $\Box $

Another interesting relation appears if we multiply (\ref{rB1})
by $x_{n-1}$
from the left and subtract  (\ref{rA2}) with $i=n-1$
multiplied from the right by $x_n^2:$
\begin{equation}
x_{n-1}x_nx_{n-1}x_n^2\equiv _n\alpha x_{n-1}x_n^2x_{n-1}x_n,
\label{Kom}
\end{equation}
in which case $\alpha =p_{n-1 n}p_{nn}^{[3]}\neq 0.$
\addtocounter{nom}{1}
\begin{lemma}
For $k<s<m\leq n$
the following relation is valid.
\begin{equation}
w_{km}w_{ks}\equiv _{k+1}\varepsilon w_{ks}w_{km},\ \ \ \varepsilon \neq 0.
\label{u51}
\end{equation}
\label{u5}
\end{lemma}
{\it Proof}. Let us use downward induction on $k.$
For this we first transpose the second letter $x_k$
of $w_{km}w_{ks}$
as far to the left as possible by means of (\ref{rA3}),
and then change the onset  $x_kx_{k+1}x_k$
according to (\ref{cong2}). We get
\begin{equation}
w_{km}w_{ks}\equiv _{k+1}\alpha x_k^2
(w_{k+1 m}w_{k+1 s}), \ \ \ \alpha \neq 0.
\label{u52}
\end{equation}

For $k+1<s$
we apply the inductive supposition to the word in the parentheses
and then by (\ref{cong2}) and (\ref{rA3}) transpose $x_k$
to its former position.

The case $k+1=s,$
the basis of the induction on $k,$
we prove by downward induction on $s.$

Let $k+1=s=n-1.$
Then $m=n.$
Let us show firstly that
\begin{equation}
\underline{x_{n-1}x_n^2x_{n-1}x_n}x_nx_{n-1}\equiv _n
\alpha x_{n-1}x_{n}^2x_{n-1}^2x_n^2+
\beta x_{n-1}x_nx_{n-1}^2x_n^3, \ \ \alpha \neq 0.
\label{soo}
\end{equation}
For this in the left hand side we transpose the first letter $x_n$
by means of (\ref{Kom}) to the penultimate
position, and then replace the ending $x_n^3x_{n-1}$
by (\ref{rB1}). We get a linear combination of three
words. One of them equals the second word of (\ref{soo}),
while two other have the following forms.
$$
x_{n-1}x_n\underline{x_{n-1}x_nx_{n-1}}x_n^2,\
\underline{x_{n-1}x_nx_{n-1}x_n^2}x_{n-1}x_n.
$$
The former word by (\ref{rA2}) transforms into the form (\ref{soo}).
The latter one, after the application of (\ref{Kom})
and the replacing of $x_{n-1}x_nx_{n-1}$
by (\ref{rA2}), will have an additional term
$\underline{x_{n-1}x_n^3}x_{n-1}^2x_n$
to which it is possible to apply (\ref{cong1}).
The direct calculation of the coefficients
shows that $\alpha =p_{n-1 n}p_{nn}\neq 0.$

Now let us multiply (\ref{soo}) by  $x_{n-2}^2$
from the left and use (\ref{rA2}) with $i=n-2.$
We get that $w_{n-2\, n}w_{n-2\, n-1}$
with respect to $\equiv _{n-1}$
equals
\begin{equation}
\gamma x_{n-2}x_{n-1}x_n^2
\underline{x_{n-2}x_{n-1}^2}x_n^2+
\delta x_{n-2}x_{n-1}x_n
\underline{x_{n-2}x_{n-1}^2}x_n^3,\ \gamma \neq 0.
\label{osn}
\end{equation}
Let us apply (\ref{cong1}) and then (\ref{cong2})
and (\ref{cong1}) to the second word. We get that this
word with respect to $\equiv _{n-1}$
equals zero. The first word after application of (\ref{rA2})
takes up the form
$$
\varepsilon w_{n-2\, n-1}w_{n-2\, n}+
\varepsilon ^{\prime } \underline{w_{n-2\, n}x_{n-1}^2}x_{n-2}x_n^2,\ \ \
\varepsilon  \neq 0.
$$
Thus, by Lemma \ref{ong4}, the  basis of the induction on
$s$ is proved.

Let us carry out the inductive step. Let $k+1=s<n-1.$
If $m>s+1=k+2$
then by the inductive supposition on $s$
we may write
$$
w_{km}w_{ks}=(w_{km}w_{k k+2})x_{k+1}
\equiv _{k+1}\alpha w_{k k+2}w_{km}x_{k+1}=
$$
\begin{equation}
\beta w_{k k+2}\underline{x_kx_{k+1}x_{k+2}x_{k+1}}w_{k+3\, m}.
\label{rav1}
\end{equation}
Taking into account (\ref{cong4})
we may neglect the words starting with $x_{k+1}^2,$ $x_{k+2}$
while  transforming the underlined part:
\begin{equation}
x_k\underline{x_{k+1}x_{k+2}x_{k+1}}
\equiv \gamma \underline{x_kx_{k+1}^2}x_{k+2}
\equiv \delta x_{k+1}x_kx_{k+1}x_{k+2}.
\label{rav2}
\end{equation}
In this way (\ref{rav1}) is transformed into (\ref{u51}).

If $m=s+1=k+2<n$
then the relation (\ref{u52}) takes up the form
$$
w_{km}w_{ks}\equiv _{k+1}\alpha x_k^2(w_{k+1 k+2}w_{k+1 k+3})
x_{k+2}x_{k+1}.
$$
Let us apply the inductive supposition with $k^{\prime }=k+1,$
$s^{\prime }=k+2,$ $m^{\prime }=k+3$
to the  word in the parentheses. We get
$$
w_{km}w_{ks}\equiv _{k+1}\alpha \varepsilon ^{-1} x_k^2
w_{k+1 k+3}w_{k+1 k+3}\underline{x_{k+2}^2x_{k+1}},
$$
or after an evident replacement
$$
w_{km}w_{ks}\equiv _{k+1}\gamma x_k^2w_{k+1 k+3}w_{k+1 k+2}\cdot
x_{k+1}x_{k+2}+\delta x_k^2 w_{k+1 k+3}^2x_{k+1}x_{k+2}^2.
$$
In both terms we may transpose one letter $x_k$
to its former position by means of (\ref{cong2}) and (\ref{rA3}). We get
\begin{equation}
w_{km}w_{ks}\equiv _{k+1}\gamma ^{\prime }
\underline{w_{k k+3}w_{k k+1}}x_{k+2}+
\delta ^{\prime }w_{k k+3}^2x_{k+1}x_{k+2}^2.
\label{u58}
\end{equation}
It is possible to apply (\ref{u51}) with
$m^{\prime }=k+3,$ $s^{\prime }=k+1$
to the first term since the case $m>s+1$
is completely considered.
Therefore it is enough to show that the second term
equals zero with respect to $\equiv _{k+1}.$
When we transpose the third letter $x_{k+1}$
as far to the left as possible we get the word
\begin{equation}
w_{k k+3}\underline{x_kx_{k+1}x_{k+2}x_{k+1}}w_{k+3\, k+3}x_{k+2}^2.
\label{u54}
\end{equation}
Taking into account (\ref{cong4})
we may neglect the words starting with $x_{k+1}$
while  transforming the underlined part:
\begin{equation}
x_k\underline{x_{k+1}x_{k+2}x_{k+1}}\equiv
x_{k+2}\underline{x_kx_{k+1}^2}\equiv
x_{k+2}x_{k+1}x_kx_{k+1}.
\label{u59}
\end{equation}
Therefore the word (\ref{u54}) equals $w_{k k+1}\underline{w_{k k+3}x_{k+2}^2}$
with respect to $\equiv _{k+1}$
and it remains only to apply  Lemma \ref{ong4} twice. \hfill $\Box $
\addtocounter{nom}{1}
\begin{lemma} The set $B$
satisfies the conditions of Lemma $\ref{tvp}.$
\label{bn7}
\end{lemma}
{\it Proof}. By Lemmas \ref{bn6} and \ref{kri1}
it is sufficient to show that in $U_P^b(\frak g)$
all words of the form $u_0, \ldots , u_6$
are linear combinations of lesser ones. The words $u_0$
are diminished by (\ref{rA3}). The words $u_1, u_2$
have been presented in this way, without using $[x_{n-1}x_n^2]=0,$
in the proof of the above theorem. The relation (\ref{cong4})
shows that $u_3\equiv _{k+1}0,$ $u_4\equiv _{k+1}0.$
Lemma  \ref{u5} with $s=m-1$
yields the necessary representation for $u_5.$

Let us prove by downward induction on $k$
that
$$
u_6\stackrel{df}{=}u_{kn}^2x_n\equiv _{k+1} \varepsilon
u_{kn}x_nu_{kn},\ \ \varepsilon \neq 0.
$$
For $k=n-1$
this equality takes up the form (\ref{Kom}). Let $k<n-1.$
Let us transpose the second letter $x_k$
of $u_{kn}^2x_n$
as far to the left as possible by means of (\ref{rA3})
and then apply  (\ref{rA1}). We get
$$
u_{kn}^2x_n\equiv _{k+1}\alpha x_k^2(u_{k+1 n}^2x_n), \ \ \ \alpha \neq 0.
$$
We may apply the inductive supposition to the
term in the parentheses and then by
(\ref{rA1}), (\ref{rA3}) transpose one of $x_k$'s
to its former position. \hfill $\Box $
\addtocounter{nom}{1}
\begin{lemma}
If $p_{11}\neq 1$
then the values of polynomials $[v]^h,$
where $[v]\in B,$ $v\neq x_i$ $h\geq 1$
are not skew primitive, in particular, they  are non-zero.
\label{blpr}
\end{lemma}
{\it Proof}. Note that for $n>2$
the sub-algebra generated by $x_2,\ldots x_n$
is defined by the Cartan matrix of the type $B_{n-1}.$
This allows us to carry out the induction on $n$
with additional supposition that the statements 1 and 2
of Theorem $B_n$
are valid for lesser values of $n.$
It is convenient formally consider the sub-algebras
$\langle x_i \rangle$
as algebras of the type $B_1.$
In this case for $n=1$
the lemma and the statements 1 and 2
are correct in the evident way.
If $v$
starts with $x_k\neq x_1$
then we may directly use the inductive supposition.
If $v=u_{1m},$
one may literally repeat the arguments of Lemma
\ref{alpr} starting at the formula (\ref{pov1}). Let $v=w_{1m}.$
If $m>2$
then by Lemma  \ref{bn1q} we have $w_{1m}=[x_1[w_{2 m}]].$
This provides a possibility to repeat the same  arguments
of Lemma \ref{alpr}  with $w$
in place of $u.$

Consider the last case $v=w_{12}.$
By Lemma \ref{bn1q} we have
\begin{equation}
[w_{12}]=[w_{13}]x_2-p(w_{13},x_2)x_2[w_{13}],
\label{pov3}
\end{equation}
\begin{equation}
[w_{13}]=x_1[w_{23}]-p(x_1,w_{23})[w_{23}]x_1.
\label{pov4}
\end{equation}
Applying the coproduct first to (\ref{pov4}) then to (\ref{pov3})
we may find the sum  $\Sigma $
of all tensors $w^{(1)}\otimes w^{(2)}$
of $\Delta ([w_{12}])$
with deg$_1(w^{(2)})=1,$  deg$_k(w^{(2)})=0,$ $k>1$
(in much the same way as (\ref{pov2})):
\begin{eqnarray}
& \Sigma =(\varepsilon g_1[w_{23}]\otimes x_1)(x_2\otimes 1)-
p(w_{13},x_2)(x_2\otimes 1)(\varepsilon g_1[w_{23}]\otimes x_1)=
\nonumber \\
& \varepsilon g_1([w_{23}]x_2-p(w_{13},x_2)p(x_2,x_1)x_2[w_{23}])\otimes x_1.
\label{bb}
\end{eqnarray}
For
$n>2,$
taking into account first
the bicharacter property  of  $p,$
then the equality
$[x_2[w_{23}]]=$ $x_2[w_{23}]-$ $p(x_2,w_{23})[w_{23}]x_2,$
and next the following  relations
$p_{ij}p_{ji}=1,$
$i-j>1;$
$p_{11}^{-1}=$
$p_{12}p_{21}=$
$p_{22}^{-1}=$
$p_{23}p_{32},$
we may write
\begin{equation}
\Sigma =\varepsilon g_1(-p(w_{13},x_2)p_{21}[x_2w_{23}]+
(1-p_{11}^{-1})[w_{23}]\cdot x_2)\otimes x_1.
\label{cb}
\end{equation}

Consider the left hand side of this tensor on applying
the inductive supposition. Note that $x_2w_{23}$
is a standard word and $[x_2w_{23}]$
equals $[x_2[w_{23}]].$
This super-letter is non-hard in $U_P(\frak g)$
since $x_2w_{23}$
contains the sub-word $x_2^2x_3.$
Thus $[x_2w_{23}]$
is a linear combination of monotonous  non-decreasing super-words
in lesser super-letters.
Among these super-words there is no  $[w_{23}]\cdot x_2$
since $x_2>x_2w_{23}.$
On the other hand, $[w_{23}]\cdot x_2$
is a monotonous non-decreasing super-word and
hence its value in $U_P(\frak g)$
is a basis element. Therefore for $n>2$
the left hand side $W$
of $\Sigma $
is non-zero.

For $n=2,$
by the definition $w_{23}=x_2,$ $w_{13}=x_1x_2,$
and the equality (\ref{bb}) takes up the form
$\Sigma =\varepsilon g_1(1-p_{12}p_{22}p_{21})x_2^2\otimes x_1.$
Since $1\neq p_{11}^{-1}=p_{12}p_{21}=p_{22}^{-2},$
we get $(1-p_{12}p_{22}p_{21})=1-p_{22}^{-1}\neq 0.$
Therefore in this case $\Sigma \neq 0$ as well.

By [{\bf \ref{pKh1}}, Corollary 10] the sub-algebra
generated by $x_2,\ldots ,x_n$
has no zero divisors. In particular $W^h\neq 0$
and $\Sigma ^h\neq 0$
in any case.

It remains to note that for $n>1$
the sum of all tensors $w^{(1)}\otimes w^{(2)}$
of $\Delta([w_{12}]^h)$
such that deg$_1(w^{(2)})=h,$  deg$_k(w^{(2)})=0,$ $k>1$
equals $\Sigma ^h,$
hence $[w_{12}]^h$
can not be skew-primitive. \hfill $\Box $

{\it Proof} of Theorem $B_n.$
Since none of $u_{km},$ $w_{km}$
contains sub-words (\ref{vor}),
Lemmas  \ref{bn7}, \ref{tver}, \ref{tvp}
imply the first statement.

If $[v]\in B$
is of finite height then by Lemma \ref{mbn} and the homogeneous  version of
Definition \ref{h1} we have $[v]^h=0.$
For $p_{11}\neq 1$
this contradicts Lemma \ref{blpr}.

Along similar lines, by Lemma \ref{prim}, every skew primitive
homogeneous element has the form $[v]^h.$
This, together with Lemma  \ref{blpr},
proves the fourth statement and, for $p_{11}\neq 1,$
the second one too.

If $p_{11}=1$
then by (\ref{sux}) we have $p_{nn}^2=1,$ $p_{ii}=1,$ $i<n.$
Besides,  $p_{ij}p_{ji}=1$
for all $i,j.$
This means that the skew commutator is a quantum operation.
Hence all elements of $B$
are skew primitive. In the case $p_{nn}=1$
these elements span a colour Lie algebra, while in the case $p_{nn}=-1$
they span a colour Lie super-algebra.
Now as in Theorem $A_n,$
we may use the  $PBW$-theorem
for  the colour Lie super-algebras.

The third statement will follow Theorem \ref{Gr}
and Lemmas \ref{min}, \ref{bn6} if we prove that all
super-letters (\ref{bGB}) are zero in $U_P({\frak g}).$
We have already proved that these super-letters are non-hard.
Therefore it remains to use the homogeneous version of Definition \ref{tv1}
and Lemma \ref{mbn}.  \hfill $\Box $

{\bf Theorem $\bf C_n$}. {\it
Suppose that ${\frak g}$
is of the type $C_n,$
and $p_{ii}\neq -1,$ $1\leq i\leq n,$ $p_{n-1 n-1}^{[3]}\neq 0.$
Denote by $B$
the set of the following super-letters:
\begin{equation}
\begin{matrix}
[u_{km}]\hfill &\stackrel{df}{=}&[x_kx_{k+1}\ldots x_m],\hfill &
1\leq k\leq m\leq n;\hfill \cr
[v_{km}]\hfill &\stackrel{df}{=}&[x_kx_{k+1}\ldots x_n\cdot x_{n-1}\ldots x_m],
\hfill &1\leq k<m<n;\hfill \cr
[v_k]\hfill &\stackrel{df}{=}&[u_{k\, n-1}u_{kn}],
\hfill &1\leq k<n.\hfill \end{matrix}
\label{CN11}
\end{equation}
The statements given below are valid.

$1.$
The values of the super-letters $(\ref{CN11})$
in $U_P({\frak g})$
form the PBW-generators set.

$2.$
Each  of these super-letters has the infinite height in $U_P({\frak g}).$

$3.$
The following relations with $(\ref{Gr3})$ form a
Groebner--Shirshov system that determines the
crystal basis of $U_P({\frak g}).$
\begin{equation}
\begin{matrix}
[u_0]\hfill &\stackrel{df}{=}&[x_kx_m]=0,\hfill &
1\leq k<m-1<n;\hfill \cr
[u_1]\hfill &\stackrel{df}{=}&
[u_{km}x_{k+1}]=0,\hfill &
1\leq k<m\leq n,\ (k,m)\neq (n-2,n);\hfill \cr
[u_2]\hfill &\stackrel{df}{=}&
[u_{km}u_{k m+1}]=0,\hfill &
1\leq k\leq m<n-1;\hfill \cr
[w_3]\hfill &\stackrel{df}{=}&
[v_{km}x_{k+1}]=0,\hfill &
1\leq k<m<n,\ k\neq m-2;\hfill \cr
[w_4]\hfill &\stackrel{df}{=}&
[v_{k k+1}x_{k+2}]=0,\hfill &
1\leq k<n-1;\hfill \cr
[w_5]\hfill &\stackrel{df}{=}&
[v_{km}v_{k m-1}]=0,\hfill &
1\leq k<m-1\leq n-1;\hfill \cr
[w_6]\hfill &\stackrel{df}{=}&
[u_{k\, n-1}^3x_n]=0,\hfill &
1\leq k<n.\hfill \end{matrix}
\label{cGB}
\end{equation}

$4.$
If $p_{11}\neq 1$
then the generators $x_i$
and their powers $x_i^t, x_i^{tl^k},$
such that $p_{ii}$
is a primitive $t$-th root of 1
together with the constants $1-g,$ $g\in G$
form a basis of ${\frak g}_P=L(U_P({\frak g})).$
Here $l$
is the characteristic of the ground field.

$5.$
If  $p_{11}=1$
then the elements $(\ref{CN11})$
and in the case of prime characteristic $l$
theirs $l^k$-th
powers, together with the constants $1-g,$ $g\in G$
form a basis of ${\frak g}_P.$}

\

In the case  $C_n$
the algebra $U_P^b(\frak g)$
is defined by the same relations (\ref{rA1}), (\ref{rA2}), (\ref{rA3}),
where in (\ref{rA2})  the last relation, $i=n-1,$
is replaced with
\begin{eqnarray}
&\ &x_{n-1}^3x_n=p_{n-1 n}p^{[3]}_{n-1 n-1}x_{n-1}^2x_nx_{n-1}+
\nonumber \\
&\ &-p_{n-1 n}^2p_{n-1 n-1}p^{[3]}_{n-1 n-1}x_{n-1}x_nx_{n-1}^2+
p_{n-1 n}^3p_{n-1 n-1}^3x_nx_{n-1}^3.
\label{rC2}
\end{eqnarray}
By Corollary  \ref{sled} we get the existence conditions
\begin{eqnarray}
& & p_{ii}=p_{11},\ p_{i-1i}p_{ii-1}=p_{11}^{-1},\ 1<i<n,\nonumber \\
& & p_{n-1n}p_{nn-1}=p_{nn}^{-1}=p_{n-1n-1}^{-2};\ p_{ij}p_{ji}=1,\ i-j>1.
\label{suxc}
\end{eqnarray}
Therefore the following relations are correct
\begin{eqnarray}
& & x_ix_{i+1}^2\equiv _{i+1}0, \ \ \ \ \ \ \ \ \ \ \ \ \  1\leq i<n;
\label{cong1c} \\
& & x_ix_{i+1}x_i\equiv _{i+1}\alpha x_i^2x_{i+1},\ \ \
1\leq i<n-1, \ \alpha \neq 0;
\label{cong2c} \\
& & x_{n-1}x_nx_{n-1}^2\equiv _n\alpha x_{n-1}^3x_n+
\beta x_{n-1}^2x_nx_{n-1},\ \  \alpha ,\ \beta \neq 0.\label{cong31c}
\end{eqnarray}
The left multiplication by $x_{n-2}$
of the last relation implies
\begin{equation}
x_{n-2}x_{n-1}x_nx_{n-1}^2\equiv _{n-1}0.
\label{ng3c}
\end{equation}
\addtocounter{nom}{1}
\begin{lemma}
The brackets in $[v_{km}], [v_k]$
are set according to the following recurrence formulae,
where by the definition $v_{k n}=u_{k n}.$
\begin{equation}
\begin{matrix}
[v_{km}]\hfill \!\!\!&=&[x_k[v_{k+1 m}]],\hfill &\hbox{if}\
1\leq k<m-1<n-1;\hfill \cr
[v_{k k+1}]\hfill \!\!\!&=&[[v_{k k+2}]x_{k+1}],\hfill &
\hbox{if}\ 1\leq k<n-1;\hfill \cr
[v_k]\hfill \!\!\!&=&[[u_{k\, n-1}][u_{kn}]],\hfill &\hbox{if}\ 1\leq k<n. \hfill \end{matrix}
\label{sc}
\end{equation}
\label{cn1}
\end{lemma}
{\it Proof}. It is enough to use the properties 6s, 1s and 2s. \hfill $\Box $
\addtocounter{nom}{1}
\begin{lemma}
If $[u], [v]\in B$
then one of the following  statements is valid.

$1)\ [[u][v]]$
is not a standard nonassociative word;

$2)\  uv$
contains a sub-word of one of the types
$u_0, u_1, u_2, w_3, w_4, w_5, w_6;$

$3)\  [[u][v]]\in B.$
\label{cn6}
\end{lemma}
{\it Proof}. The first two formulae (\ref{sc})
coincide with (\ref{sb}) up to replacement of $v$
with $w$
provided $k+1\neq n>m.$
Obviously for $m<n$
the inequality $v_{km}>v_{rs}$
is equivalent to $w_{km}>w_{rs},$
while  $v_{km}>u_{rs}$
is equivalent to $w_{km}>w_{rs}.$
Hence Lemmas (\ref{bn2}), (\ref{bn3}), (\ref{bn4})
are still valid under the replacement of $w$
with $v:$
\begin{equation}
\begin{matrix}
[[v_{km}][v_{rs}]]& \hbox{ is standard }\Leftrightarrow \hfill &
s\geq m>k+1=r\vee (s<m\& r=k);\hfill \cr
[[u_{km}][v_{rs}]]& \hbox{ is standard }\Leftrightarrow \hfill &
k=r\vee k=m<r;\hfill \cr
[[v_{km}][u_{rs}]]& \hbox{ is standard }\Leftrightarrow \hfill &
r=k+1<m\vee r=k+1=m=s.\hfill \end{matrix}
\label{cn}
\end{equation}
Further, $v_k>v_r$
if and only if $k<r,$
and under this condition $[[v_k][v_r]]$
is not standard since $u_{kn}>u_{r n-1}u_{rn}.$

In a similar manner $v_k>u_{rm}$
is equivalent to $k<r,$
while $v_k>v_{rm}$
is equivalent to $k\leq r.$
Therefore none of the words  $[[v_k][u_{rm}]],$ $[[v_k][v_{rm}]]$
is standard since $u_{kn}>u_{rm}$
and $u_{kn}>v_{rm},$
respectively.

For the remaining two cases we have only two possibilities
\begin{equation}
\begin{matrix}
[[u_{km}][v_r]] & \hbox{ is standard }\Leftrightarrow \hfill &
r=k\leq m<n;\hfill \cr
[[v_{km}][v_r]] & \hbox{ is standard }\Leftrightarrow \hfill &
r=k+1\& k<m-1.\hfill \end{matrix}
\label{cn1q}
\end{equation}
\noindent The treatment in turn of the eight possibilities
(\ref{cn}), (\ref{cn1q}) proves the lemma. \hfill $\Box $
\addtocounter{nom}{1}
\begin{lemma}
If a super-word $W$
equals one of the super-letters $(\ref{cGB})$
or $[v]^h,$ $[v]\in B,$ $h\geq 1,$
then its constitution does not equal the constitution
of any word in less then $W$
super-letters from $B.$
\label{mcn}
\end{lemma}
{\it Proof}. The proof is akin to Lemma \ref{man}
with the following tableaux:
\begin{equation}
\begin{matrix}
[u_{km}]^h,\ [u_{km}x_{k+1}], \ [u_{km}u_{km+1}] \hfill &
\hbox{deg}_k(u)\leq \hbox{deg}_{m+1}(u); \hfill \cr
[v_{km}]^h,\ [v_{km}x_{k+1}], \ [v_{km}v_{km-1}] \hfill &
2\hbox{deg}_k(u)\leq \hbox{deg}_{m-1}(u); \hfill \cr
[v_{kk+1}x_{k+2}] \hfill &
\hbox{deg}_k(u)=0; \hfill \cr
[v_k]^h \hfill &
\hbox{deg}_k(u)\leq \hbox{deg}_n(u); \hfill \cr
[u_{kn-1}^3x_n]\hfill &
\hbox{deg}_k(u)\leq 2\hbox{deg}_n(u).\hfill \end{matrix}
\label{mcf}
\end{equation}

\hfill $\Box $
\addtocounter{nom}{1}
\begin{lemma}
If $y=x_i,\ m-1\neq i>k$
or $y=x_i^2,\ m-1=i>k$
then
\begin{equation}
v_{km}y\equiv _{k+1}0.
\label{cong4c}
\end{equation}
\label{ong4c}
\end{lemma}
{\it Proof}. For $i<m-1,$
we may transpose $y$
by means of (\ref{rA3}) to the left across $x_n^2$
and then use Lemma \ref{ong4a} with $m^{\prime }=n-1.$

If $y=x_i^2,\ $ $m-1=i>k$
then by the above case, $i<m-1,$
we get
\begin{equation}
v_{km}y=v_{k m+1}\underline{x_mx_{m-1}^2}=
\underline{v_{k m+1}x_{m-1}}
(\alpha x_mx_{m-1}+\beta x_{m-1}x_m)\equiv _{k+1}0,
\label{ssc}
\end{equation}
where by definition $v_{k n}=u_{kn}$
and $u_{kn}x_{n-2}\equiv _{n-2}0,$
while $n-2=i>k$.

If $y=x_i,$ $i=m>k$
then for $m=n-1$
we may use the inequality (\ref{ng3c}), while for $m<n-1$
we have $v_{km}y=v_{k m+1}y_1$
where $y_1=x_m^2.$
Hence we may use (\ref{ssc}) replacing $m$
by $m+1.$

If $y=x_i,\ i>m>k$
then by (\ref{rA3}) we get
$v_{km}y=\alpha v_{k i+1}\underline{x_i x_{i-1}x_i}\cdot w.$
Changing the underlined by (\ref{rA1}), we may apply the
previously considered  cases:  $m^{\prime }-1=i^{\prime },$
where $m^{\prime }=i+1,$ $i^{\prime }=i;$
and $i^{\prime }<m^{\prime }-1,$
where $m^{\prime }=i+1,$ $i^{\prime }=i-1.$   \hfill $\Box $

\

If we multiply (\ref{rC2}) by $x_n$
from the right and subtract  (\ref{rA1}) with $i=n-1$
multiplied from the left by $x_{n-1}^2,$
then by means of $p_{n-1 n-1}^{-2}=p_{n n-1}p_{n-1 n}=p_{n n}^{-1}$
we get
\begin{equation}
x_{n-1}^2\underline{x_nx_{n-1}x_n}\equiv _n p_{n-1 n}(p_{n-1 n-1}^{[3]}
x_{n-1}x_nx_{n-1}^2x_n-
p_{n-1 n-1}x_{n-1}^2x_n^2x_{n-1}).
\label{Koc}
\end{equation}
Let us first multiply this relation by $x_{n-2}^2$
from the left and then apply (\ref{rA1}) to the
underlined sub-word. Taking into account   the relation
$x_{n-2}^2x_{n-1}^3\equiv _{n-1}0,$
we get that the left hand side of the multiplied (\ref{Koc}) equals
$p_{n-1 n}p_{n n}(1+p_{n n})^{-1}x_{n-2}^2x_{n-1}^2x_n^2x_{n-1}$
up to  $\equiv _{n-1},$
i.e. it is proportional to the second term of the right
hand side. As a result the relation below with
$\alpha =p_{n-1 n-1}^{-1}(1+p_{n n})\neq 0$
is correct.
\begin{equation}
x_{n-2}^2x_{n-1}^2x_n^2x_{n-1}\equiv _{n-1}
\alpha x_{n-2}^2x_{n-1}x_nx_{n-1}^2x_n.
\label{Komc}
\end{equation}
\addtocounter{nom}{1}
\begin{lemma}
If $k<s<m\leq n$
and as above $v_{kn}=u_{kn}$
then
\begin{equation}
v_{km}v_{ks}\equiv _{k+1}\varepsilon v_{ks}v_{km},\ \ \varepsilon \neq 0,
\label{c51}
\end{equation}
\label{c5}
\end{lemma}
{\it Proof}.
Let us use downward induction on $k.$
For this we first  transpose the second letter $x_k$
of $v_{km}v_{ks}$
as far to the left as possible by means of (\ref{rA3}),
and then change the onset $x_kx_{k+1}x_k$
according to (\ref{cong2c}). We get
\begin{equation}
v_{km}v_{ks}\equiv _{k+1}\alpha x_k^2
(v_{k+1 m}v_{k+1 s}), \ \ \ \alpha \neq 0.
\label{u52c}
\end{equation}

For $k+1<s$
we may apply the inductive supposition to the word
in the parentheses, and then transpose $x_k$
to its former position by (\ref{cong2c}), (\ref{rA3}).

For $k+1=s$
we will use downward induction on $s.$

Let $k+1=s=n-1.$
In this case $m=n$
and (\ref{u52c}) becomes:
$$
v_{n-2\, n}v_{n-2\, n-1}\equiv _{n-1}\beta x_{n-2}^2
(x_{n-1}\underline{x_nx_{n-1}x_n}x_{n-1}).
$$
Let us replace the underlined part according to (\ref{rA1}).
Since $x_{n-2}^2x_{n-1}x_n^2\equiv _n0,$
we may continue by (\ref{Komc}):
$$
\equiv _{n-1}\beta _1x_{n-2}^2x_{n-1}^2x_n^2x_{n-1}\equiv_{n-1}
\beta _2\underline{x_{n-2}^2x_{n-1}}x_nx_{n-1}^2x_n\equiv _{n-1}
$$
$$
\beta _3x_{n-2}x_{n-1}\underline{x_{n-2}x_n}x_{n-1}^2x_n\equiv _{n-1}
\beta _4x_{n-2}x_{n-1}x_n\underline{x_{n-2}x_{n-1}^2}x_n.
$$
With the help of (\ref{rA1}) we get
$$
=\varepsilon v_{n-2\, n-1}v_{n-2\, n}+\beta _5 x_{n-2}
\underline{x_{n-1}x_nx_{n-1}^2}x_{n-2}x_n,\ \ \varepsilon \neq 0.
$$
By (\ref{cong31c}) and
(\ref{cong1c}) we see that the second term equals zero
up to $\equiv _{n-1}.$

The inductive step on $s$
coincides the inductive step on $s$
in Lemma \ref{u5} up to replacing
both the citations of Lemma \ref{ong4}
with the citations of Lemma \ref{ong4c} and $w$
with $v.$  \hfill $\Box $
\addtocounter{nom}{1}
\begin{lemma}
The set  $B$
satisfies the Lemma $\ref{tvp}$ conditions.
\label{cn7}
\end{lemma}
{\it Proof}. According to the Super-letter Crystallisation Lemma
and Lemma \ref{bn6} it is sufficient to show that words of the form
$u_0,u_1,u_2,w_3,w_4,w_5,w_6$
are linear combinations of lesser words in $U_P(\frak g).$
The words $u_0$
are diminished by (\ref{rA3}).
The words $u_1, u_2$
have been diminished in Theorem $A_n$
since in the case $C_n$ the words $u_2$
are independent of $x_n,$ while $u_1$
depends  on $x_n$
only if $u_1=x_{n-1}x_n^2.$
The relation (\ref{cong4c}) shows  that
$w_3\equiv _{k+1}0,$  $w_4\equiv _{k+1}0.$
Lemma  \ref{c5} with $s=m-1$
gives the required representation for $u_5.$

Consider  the words $w_6.$
For   $k=n-1$
the relation (\ref{rC2}) defines the required  decomposition. Let $k<n-1.$
Since $x_1,\ldots , x_{n-1}$
generate a sub-algebra
of the type $A_{n-1},$
the crystal decomposition of $u_{kn-2}^3x_{n-1}$
has the form
\begin{equation}
u_{kn-2}^3x_{n-1}=\sum \alpha u_{m_1s_1}u_{m_2s_2}\cdots u_{m_ts_t},
\label{cry}
\end{equation}
where $u_{m_1s_1}\leq u_{m_2s_2}\leq \ldots \leq u_{m_ts_t},$
that is $m_1\geq m_2\geq \ldots \geq m_t$
and  $s_i\geq s_{i+1}$
if $m_i=m_{i+1}.$
In particular, if $m_1=k$ then $m_2=\ldots =m_t=k$
and, due to the homogeneity, $t=3,$ $s_1=n-1,$ $s_2=s_3=n-2.$
Therefore
\begin{equation}
u_{kn-2}^3x_{n-1}\equiv _{k+1}\varepsilon u_{kn-1}u_{kn-2}^2.
\label{1eq}
\end{equation}
Along similar lines, the following relations are valid as well
\begin{equation}
u_{kn-2}^3x_{n-1}^2\equiv _{k+1}\mu u_{kn-1}^2u_{kn-2},
\ \ \ u_{k\, n-2}^2x_{n-1}^3\equiv _{k+1}0.
\label{2eq}
\end{equation}

Now let us multiply  (\ref{rA1}) with $i=n-2$ by $x_{n-1}$
from the right, and then add to the result the same relation
multiplied by  $p_{n-2\, n-1}(1+p_{n-1 n-1})x_{n-1}$
from the left. We get the following relation
with $\alpha =p_{n-2\, n-1}^2p_{n-1n-1}^{[3]}\neq 0.$
\begin{equation}
x_{n-2}x_{n-1}^3=
\alpha x_{n-1}^2x_{n-2}x_{n-1}+\beta x_{n-1}^3x_{n-2},
\label{3c}
\end{equation}

Further, we may write
\begin{equation}
u_{k\, n-1}^3=\beta _1u_{k\, n-2}u_{k\, n-3}
\underline{x_{n-1}x_{n-2}x_{n-1}}u_{k\, n-1},\ \ \beta _1\neq 0,
\label{3ca}
\end{equation}
where for $k=n-2$
the term  $u_{k\, n-3}$
is absent. Let us apply (\ref{rA1}) with $i=n-2$
to the underlined word. Since $u_{k\, n-2}u_{k\, n-3}x_{n-1}^2\equiv _{n-1}0,$
we have got
\begin{equation}
u_{k\, n-1}^3\equiv _{n-1}
\beta _2u_{k\, n-2}^2u_{k\, n-3}\underline{x_{n-1}^2x_{n-2}x_{n-1}}.
\label{3cb}
\end{equation}
Let us apply (\ref{3c}). Taking into account
the second of (\ref{2eq}) we get
\begin{equation}
u_{k\, n-1}^3\equiv _{k+1}
\beta _3u_{k\, n-2}^3x_{n-1}^3.
\label{3cd}
\end{equation}
Let us multiply this relation from the right by  $x_n.$
By (\ref{rC2}) we have
\begin{equation}
u_{k\, n-1}^3x_n\equiv _{k+1}
\alpha \underline{u_{k\, n-2}^3x_{n-1}}x_nx_{n-1}^2+
\beta \underline{u_{k\, n-2}^3x_{n-1}^2}x_nx_{n-1}.
\label{3ce}
\end{equation}
By means of (\ref{1eq}) and (\ref{2eq}) we have got
$$
u_{k\, n-1}^3x_n\equiv _{k+1}
\alpha _1u_{k\, n-1}x_nu_{k\, n-2}^2x_{n-1}^2+
\beta_1u_{k\, n-1}^2x_nu_{k\, n-2}x_{n-1},
$$
and both of these words are less than $u_{k\, n-1}^3x_n.$ \hfill $\Box $
\addtocounter{nom}{1}
\begin{lemma}
If  $p_{11}\neq 1$
then the values of $[v]^h,$
where $[v]\in B,$ $v\neq x_i,$ $h\geq 1$
are not skew primitive. In particular they are non-zero.
\label{clpr}
\end{lemma}
{\it Proof}. Note that for $n>3$
the algebra generated by $x_2,\ldots x_n$
is a sub-algebra of the type $C_{n-1}.$
Therefore we may use induction on $n$
with additional supposition that the theorem
statements 1 and 2 are valid for the lesser  values of $n.$
We will formally consider the sub-algebra generated by $x_{n-1}, x_n$
as an algebra of the type $C_2,$
and the sub-algebra generated by $x_n$
as an algebra of type $C_1.$
In this case for $n=1$
the present lemma and the statements 1 and 2 are valid in obvious way.

If the first letter $x_k$
of $v$
is less than $x_1$
then we may use the inductive supposition directly.
If $v=u_{1m}$
then one may literally repeat arguments of
Lemma \ref{alpr} starting at (\ref{pov1}).

If $v=v_{1m}$
and $n>3$
then we may  repeat arguments of Lemma
\ref{blpr} starting at (\ref{pov3}) up to replacing $w$
with $v.$
For $n=3$
in these arguments the formula (\ref{cb}) assumes  the form
\begin{equation}
\Sigma =\varepsilon g_1(-p(v_{13},x_2)p_{21}[x_2^2x_3]+
(1-p_{11}^{-1})[x_2x_3]\cdot x_2)\otimes x_1.
\label{cbc}
\end{equation}
Therefore the left component of the tensor $\Sigma $
is a non-zero linear combination of  the basis elements.
For  $n=2$
the set  $B$
has no elements $v_{1m}$
at all.

Consider the last case,  $v=v_1=[u_{1n-1}^2x_n].$
Let $S_k$
be the sum of all tensors of $\Delta ([u_{kn}])=$
$\sum u^{(1)}\otimes u^{(2)}$
with deg$_n(w^{(1)})=1,$ deg$_k(w^{(1)})=0,$ $k<n.$
Evidently $S_n=x_n\otimes 1.$
Let us show by downward induction on $k$ that
$S_k=(1-p_{11}^{-1})g(u_{kn-1})x_n\otimes [u_{kn-1}]$
at $k<n.$
We have
\begin{equation}
\Delta ([u_{kn}])=\Delta (x_k)\Delta ([u_{k+1n}])-p(x_k, u_{k+1n})
\Delta ([u_{k+1n}])\Delta (x_k).
\label{cc1}
\end{equation}
Consequently,
\begin{equation}
S_k=(g_k\otimes x_k) S_{k+1}-p(x_k, u_{k+1n})S_{k+1}(g_k\otimes x_k).
\label{cc2}
\end{equation}
This implies the required formula since
by  (\ref{suxc})  at  $k<n-1$
we have
$$p(x_k,u_{k+1n})p(x_n,x_k)=p(x_k,u_{k+1n-1}),$$
while at  $k=n-1$
we have $p(x_{n-1},x_n)p(x_n,x_{n-1})=p_{11}^{-1}.$

In a similar manner, consider the sum $S$
of all tensors  of
$\Delta ([u_{kn}^2x_n])=$
$\sum w^{(1)}\otimes w^{(2)}$
with deg$_n(w^{(1)})=1,$ deg$_i(w^{(1)})=0,$ at $i<n.$
\begin{equation}
\Delta ([[u_{1n-1}][u_{1n}]])=
\Delta ([u_{1n-1}])\Delta ([u_{1n}])-p(u_{1n-1}, u_{1n})
\Delta ([u_{1n}])\Delta ([u_{1n-1}]).
\label{cc3}
\end{equation}
Since we now  $S_1,$  we may calculate  $S:$
\begin{eqnarray}
& & S=(g(u_{1n-1})\otimes [u_{1n-1}])S_1-p(u_{1n-1}, u_{1n})S_1
(g(u_{1n-1})\otimes [u_{1n-1}])=
\nonumber \\
& & (1-p_{11}^{-1})g(u_{1n-1}^2)x_n\otimes
(1-p(u_{1n-1}, u_{1n})p(x_n, u_{1n-1}))[u_{1n-1}]^2.
\label{cc4}
\end{eqnarray}
By (\ref{suxc}), using the bicharacter property of $p,$ we have
$$1-p(u_{1n-1}, u_{1n})p(x_n, u_{1n-1})=
1-p(u_{1n-1}, u_{1n-1})p_{n-1n}p_{nn-1}=
$$
$$
1-p_{n-1n-1}p_{n-1n-1}^{-2}=
1-p_{11}^{-1}\neq 0.$$
Because of this, $S\neq 0$
and the sum of all tensors $w^{(1)}\otimes w^{(2)}$
with deg$_n(w^{(1)})=h,$ deg$_k(w^{(1)})=0,$ $k<n$
of the basis decomposition of  $\Delta([v_1]^h)$
equals $S^h\neq 0.$ Therefore $[v_{1}]^h$
is not skew primitive. \hfill $\Box $

\

{\it Proof} of Theorem  $C_n.$
For the first statement it will suffice to prove that all
super-letters (\ref{CN11}) are hard in $U_P({\frak g}).$
Since none of  $u_{km},$ $v_{km}$
contains a sub-word  (\ref{vor}), Lemma  \ref{tver}
implies that  $[u_{km}],$ $[v_{km}]$
are hard.

If $[v_k]$
is not hard then, by the homogeneous version of Definition
\ref{tv1},  its value is a polynomial in lesser hard super-letters.
In line with Lemmas \ref{cn7} and  \ref{tvp},
all hard super-letters belong to $B.$
Therefore, by Lemma \ref{mcn}, $[v_k]=0.$
Since deg$_n(v_k)=1$ and  deg$_{n-1}(v_k)=2,$
the equality  $[v_k]=0$
is valid in the algebra $C^{\prime }$
which  is defined by all relations  of $U_P({\frak g}),$
but ones of degree greater than 1 in $x_n$
and ones of degree greater than 2 in $x_{n-1},$
that is in the algebra defined by  (\ref{rA1}), (\ref{rA2})
with $i<n-1,$
and (\ref{rA3}). These relations do not reverse
the order of  $x_{n-1}$
and  $x_n$
in monomials since none of them has both $x_{n-1}$ and $x_n.$
This implies that the sum of all monomials of
$[v_k]=[u_{kn-1}]\cdot [u_{kn}]-$
$p(u_{kn-1},u_{kn})[u_{kn}]\cdot [u_{kn-1}]$
in which  $x_n$
is prefixed to $x_{n-1}$
equals zero in $C^{\prime },$
that is  $[u_{kn}]\cdot [u_{kn-1}]=0.$
Especially, this equality is valid in $U_P({\frak g}).$
Since, by Theorem \ref{BW}, the
super-word $[u_{kn}]\cdot [u_{kn-1}]$
is a basis element, the first statement is proved.

If ${[v]}\in B$
is of finite height then, by Lemma \ref{mcn}
and the homogeneous version of Definition \ref{h1}, we have $[v]^h=0.$
For $p_{11}\neq 1$
this contradicts Lemma \ref{clpr}.
 In a similar manner, according to Lemma
 \ref{prim}, every skew primitive homogeneous
element has the form  $[v]^h.$
This, together with Lemma \ref{clpr},
proves the fourth statement and, for  $p_{11}\neq 1,$
 the second one too.
If $p_{11}=1$
then according to (\ref{suxc}) we have $p_{ii}=p_{ij}p_{ji}=1$
at all $i,j.$
In particular, the skew commutator is a quantum
operation. Hence all elements of $B$
are skew primitive. These elements span a colour
Lie algebra. Now, as in Theorem $A_n,$
we may use the coloured PBW theorem.

The third statement will follow from Theorem
\ref{Gr} and Lemmas \ref{min}, \ref{cn6} provided we note that all
super-letters  (\ref{cGB}) are zero in $U_P({\frak g}).$
We have proved already that these super-letters are non-hard.
So it  remains to use first the homogeneous version of Definition \ref{tv1}
and  then Lemma \ref{mdn}. \hfill $\Box $

{\bf Theorem $\bf D_n$}. {\it
Let ${\frak g}$
be of the type $D_n,$
and $p_{ii}\neq -1,$ $1\leq i\leq n.$
Denote by $B$
the set of the following super-letters:
\begin{equation}
\begin{matrix}
[u_{km}]\hfill &\stackrel{df}{=}&[x_kx_{k+1}\ldots x_m],\hfill &
1\leq k\leq m<n;\hfill \cr
[e_{km}]\hfill &\stackrel{df}{=}&[x_kx_{k+1}\ldots x_{n-2}\cdot
x_nx_{n-1}\ldots x_m],
\hfill &1\leq k<m\leq n,\hfill \cr
[e_{n-1n}]\hfill &\stackrel{df}{=}&x_n.\hfill \end{matrix}
\label{DN1}
\end{equation}
The statements given below are valid.

$1.$
The values of $(\ref{DN1})$ in $U_P({\frak g})$
form the PBW-generators set.

$2.$
Each of  the super-letters  $(\ref{DN1})$ has
infinite height in  $U_P({\frak g}).$

$3.$
The relations $(\ref{Gr3})$ together with the following
ones form a Groebner--Shirshov system that determines
the crystal basis of  $U_P({\frak g}).$
\begin{equation}
\begin{matrix}
[u_0]\hfill &\stackrel{df}{=}&[x_kx_m]=0,\hfill &
1\leq k<m-1<n,\ (k,m)\neq (n-2,n);\hfill \cr
[u_1]\hfill &\stackrel{df}{=}&
[u_{km}x_{k+1}]=0,\hfill &
1\leq k<m<n;\hfill \cr
[u_1^{\prime }]\hfill &\stackrel{df}{=}&
[x_{n-2}x_n^2]=0,\hfill & \cr
[u_2]\hfill &\stackrel{df}{=}&
[u_{km}u_{k\, m+1}]=0,\hfill &
1\leq k\leq m<n-1;\hfill \cr
[v_3]\hfill &\stackrel{df}{=}&
[e_{km}x_{k+1}]=0,\hfill &
1\leq k<m\leq n,\ n-1\neq k\neq m-2;\hfill \cr
[v_4]\hfill &\stackrel{df}{=}&
[e_{k k+1}x_{k+2}]=0,\hfill &
1\leq k<n-2;\hfill \cr
[v_4^{\prime }]\hfill &\stackrel{df}{=}&
[e_{n-3\,  n-2}x_n]=0,\hfill & \cr
[v_5]\hfill &\stackrel{df}{=}&
[e_{km}e_{k\, m-1}]=0,\hfill &
1\leq k<m-1\leq n-1;\hfill \cr
[v_6]\hfill &\stackrel{df}{=}&
[u_{km}e_{kn}]=0,\hfill &
1\leq k\leq m<n,\ n-2\leq m.\hfill \end{matrix}
\label{dGB}
\end{equation}

$4.$
If $p_{11}\neq 1,$
then the generators  $x_i,$
their powers  $x_i^t, x_i^{tl^k},$
such that $p_{ii}$
is a primitive $t$-th root of 1, together with the constants $1-g,$ $g\in G$
form a basis of  ${\frak g}_P=L(U_P({\frak g})).$
Here $l=$ char$({\bf k}).$

$5.$
If  $p_{11}=1,$
then the elements of  $B$
and, for $l>0,$ their $l^k$-th powers
together with the constants $1-g,$ $g\in G$
form a basis of  ${\frak g}_P.$}

\

In the case $D_n$
the algebra $U_P^b({\frak g})$
can be defined by the condition that
the sub-algebras $U_{n-1}$
and $U_n$
generated, respectively, by $x_1,\ldots ,x_{n-1}$
and $x_1,\ldots ,x_{n-2}, x_{n-1}^{\prime }=x_n$
are quantum universal enveloping algebras of the type $A_{n-1},$
and by the only additional relation
\begin{equation}
[x_{n-1}x_n]=0.
\label{rD3}
\end{equation}
The existence conditions take up the form
\begin{eqnarray}
& & p_{ii}=p_{nn}=p_{11},\ p_{i+1i}p_{ii+1}=
p_{n-2n}p_{nn-2}=p_{11}^{-1},\ \hbox{ if }1\leq i<n,\nonumber \\
& & p_{n-1n}p_{nn-1}=p_{ij}p_{ji}=1,\
\hbox{ if } i-j>1\& (i,j)\neq (n,n-2).
\label{suxd}
\end{eqnarray}
\addtocounter{nom}{1}
\begin{lemma}
The brackets in $(\ref{DN1})$ are set up by
the recurrence formulae
\begin{equation}
\begin{matrix}
[e_{km}]\hfill \!\!\!&=&[x_k[e_{k+1 m}]],\hfill &
\hbox{if}\ 1\leq k<m-1<n,\ k\neq n-1;\hfill \cr
[e_{k k+1}]\hfill \!\!\!&=&[[e_{k k+2}]x_{k+1}],\hfill &
\hbox{if}\ 1\leq k<n-1.\hfill \end{matrix}
\label{sd}
\end{equation}
\label{dn1q}
\end{lemma}
{\it Proof}. It is enough to use the properties 6s, 1s,  and 2s. \hfill $\Box $
\addtocounter{nom}{1}
\begin{lemma}
If $[u], [v]\in B,$
then one of the statements below is correct.

$1)\ [[u][v]]$
is not a standard nonassociative word;

$2)\  uv$
contains a sub-word of one of the types
$u_0, u_1, u_1^{\prime }u_2, v_3, v_4, v_4^{\prime }, v_5, v_6;$

$3)\  [[u][v]]\in B.$
\label{dn6}
\end{lemma}
{\it Proof}. The formulae  (\ref{sd}) coincides with
 (\ref{sb}) at  $k\neq n-1$
up to replacing  $e$
by $w.$
The inequality  $e_{km}>e_{rs}$
is set up by the same conditions,  $k<r\vee (k=r\& m<s),$
as the inequality $w_{km}>w_{rs}$
does. Likewise $u_{km}>e_{rs}$
is set up by the same condition, $k\leq r,$
as $u_{km}>w_{rs}$
does. Therefore Lemmas  \ref{bn2}, \ref{bn3}, \ref{bn4}
remain valid with $e$
in place of $w:$
\begin{equation}
\begin{matrix}
[[e_{km}][e_{rs}]]\hfill  & \hbox{ is standard }\Leftrightarrow \hfill &
s\geq m>k+1=r\vee (s<m\& r=k);\hfill \cr
[[u_{km}][e_{rs}]]\hfill  & \hbox{ is standard }\Leftrightarrow \hfill &
k=r\vee k=m<r;\hfill \cr
[[e_{km}][u_{rs}]]\hfill  & \hbox{ is standard }\Leftrightarrow \hfill &
r=k+1<m\vee r=k+1=m=s.\hfill \end{matrix}
\label{dn}
\end{equation}
By looking over all of these possibilities we get the lemma statement. \hfill $\Box $
\addtocounter{nom}{1}
\begin{lemma}
If a super-word $W$
equals one of the super-letters $(\ref{dGB})$
or $[v]^h,$ $[v]\in B,$ $h\geq 1$
then its constitution does not equal the constitution
of any super-word in less than $W$
super-letters from $B.$
\label{mdn}
\end{lemma}
{\it Proof}. The proof is similar to the one of Lemma \ref{man} with the tableaux
\begin{eqnarray}
& & \begin{matrix}
[u_{km}]^h,\hfill & [u_{km}x_{k+1}], \hfill & [u_{km}u_{k\, m+1}] \hfill &  &
\hbox{deg}_k(u)\leq \hbox{deg}_{m+1}(u); \hfill \cr
[e_{km}]^h,\hfill & [e_{km}x_{k+1}], \hfill &
[e_{km}e_{k\, m-1}],\hfill & m<n \hfill &
2\hbox{deg}_k(u)\leq \hbox{deg}_{m-1}(u); \hfill \cr
[e_{kn}]^h,\hfill & [e_{kn}x_{k+1}], \hfill & [e_{kn}e_{k\, n-1}]\hfill & &
\hbox{deg}_k(u)\leq \hbox{deg}_{m-1}(u); \hfill \end{matrix} \nonumber \\
& & \begin{matrix}
[e_{kk+1}x_{k+2}] \hfill &
\hbox{deg}_k(u)=0; \hfill  \cr
[e_{n-3n-2}x_n] \hfill &
\hbox{deg}_{n-3}(u)=0; \hfill \cr
[u_{k\, n-2}e_{kn}]\hfill &
\hbox{deg}_k(u)\leq \hbox{deg}_{n-1}(u)+\hbox{deg}_n(u);\hfill  \cr
[u_{k\, n-1}e_{kn}]\hfill &
\hbox{deg}_k(u)\leq \hbox{deg}_n(u). \hfill \end{matrix}
\label{mdf}
\end{eqnarray}
\addtocounter{nom}{1}
\begin{lemma}
If $y=x_i,$ $\ m-1\neq i>k$
or $y=x_i^2,$ $\ m-1=i>k$
then
\begin{equation}
e_{km}y\equiv _{k+1}0.
\label{cong4d}
\end{equation}
\label{ong4d}
\end{lemma}
{\it Proof}. If $i<m-1,$ $m\neq n,$
or $m=n,$ $i<n-2,$
then with the help of  (\ref{rA3})  and (\ref{rD3})
it is possible to permute $y$
to the left beyond $x_n$
and then to use Lemma  \ref{ong4a} for $U_{n-1}.$

If $m=n,$ $i=n-2$
then we may use Lemma  \ref{ong4a} for $U_n.$

If $y=x_i^2,$ $m-1=i>k$
then for $m<n$
by the above case we get
\begin{equation}
e_{km}y=e_{k\, m+1}x_mx_{m-1}^2=\underline{e_{k\, m+1}x_{m-1}}
(\alpha x_mx_{m-1}+\beta x_{m-1}x_m)\equiv _{k+1}0.
\label{ssd}
\end{equation}
For  $m=n$
we have  $e_{kn}x_{n-1}^2=$
$\alpha \underline{u_{k\, n-2}x_{n-1}^2}x_n\equiv _{n-1}0$
since the underlined part belongs to $U_{n-1}.$

If $y=x_i,$ $i=m>k$
then for $m=n$
we may use Lemma  \ref{ong4a} applied to $U_n;$
for $m=n-1$
we may use the same lemma  applied to $U_{n-1}$
provided that beforehand  we permute $x_n$
with $y$
by  (\ref{rD3}); for $m<n-1$
we may first  rewrite $e_{km}y=e_{k\, m+1}y_1,$
where $y_1=x_m^2,$
and then use   (\ref{ssd})  with  $m+1$
in place of $m.$

If  $y=x_i,$ $i>m>k$
then for $i<n$
we have $e_{km}y=\alpha e_{ki+1}\underline{x_ix_{i-1}x_i}\cdot v.$
Replacing the underlined word by  (\ref{rA1}) in  $U_{n-1},$
we may use the previously considered cases: $m^{\prime }-1=i^{\prime },$
where $m^{\prime }=i+1,$ $i^{\prime }=i;$
and  $i^{\prime }<m^{\prime }-1,$
where $m^{\prime }=i+1,$ $i^{\prime }=i-1.$
For $i=n,$
and $m=n-1$
we have $e_{k\, n-1}x_n=\alpha \underline{u_{k\, n-2}x_n^2}x_{n-1}$
and one may apply  Lemma \ref{ong4a} to $U_n.$
Finally, for $i=n$
and $m<n-1$
we get
$$
e_{km}x_n=\beta _1u_{k\, n-2}x_nx_{n-1}x_{n-2}x_n\cdot v=
\beta _2u_{k\, n-2}x_{n-1}\underline{x_nx_{n-1}x_n}\cdot v=
$$
$$
\beta _3\underline{u_{k\, n-2}x_{n-1}x_{n-2}}x_n^2\cdot v+
\beta _4u_{k\, n-2}\underline{x_{n-1}x_n^2}x_{n-2}\cdot v.
$$
One may apply first Lemma  \ref{ong4a} for  $U_{n-1}$
to the underlined sub-word of the first term, and then, after (\ref{rD3}),
Lemma  \ref{ong4a} for  $U_n$ to the second term. \hfill $\Box $
\addtocounter{nom}{1}
\begin{lemma}
If $k<s<m\leq n$
then
$e_{km}e_{ks}\equiv _{k+1}\varepsilon e_{ks}e_{km},$ $\varepsilon \neq 0.$
\label{d5}
\end{lemma}
{\it Proof}. Let us carry out downward induction on  $k.$
The largest value of $k$
equals $n-2.$
In this case  $s=n-1,$ $m=n$
and we have
\begin{eqnarray}
\underline{x_{n-2}x_n\cdot x_{n-2}}x_nx_{n-1}\equiv _n
x_{n-2}^2\underline{x_n^2x_{n-1}}=
\alpha \underline{x_{n-2}^2x_{n-1}}x_n^2\equiv _{n-1} \nonumber \\
\beta x_{n-2}x_{n-1}\underline{x_{n-2}x_n^2}\equiv _n
\varepsilon x_{n-1}x_n\cdot x_{n-1}x_{n-2}x_n.
\label{dosn}
\end{eqnarray}
Let us first transpose the second letter $x_k$
of $e_{km}e_{ks}$
as far to the left as possible by (\ref{rA3}),
and then replace the onset  $x_kx_{k+1}x_k$
 by (\ref{cong1a}). We get
\begin{equation}
e_{km}e_{ks}\equiv _{k+1}\alpha x_k^2
(e_{k+1 m}e_{k+1 s}), \ \ \ \alpha \neq 0.
\label{u52d}
\end{equation}

For  $k+1<s$
it suffices to apply the inductive supposition
to the word in the parentheses  and then by
(\ref{cong1a}) and (\ref{rA3}) to put $x_k$
to the proper place.

For $k+1=s$
one may use downward induction on $s.$
The basis of this induction, $s=n-1,$
has been proved, see  (\ref{dosn}). For $k<n-3$
the inductive step on $s$
 coincides with the one of Lemma  \ref{u5}
with  $e$
in place of $w$
since in this case the active variables $x_k,$ $x_{k+1}$
$q$-commute with $x_n.$
If $k=n-3$
then in consideration of Lemma  \ref{u5}
the variable $x_{k+1}=x_{n-2}$
is transposed across $x_n$
twice: in  (\ref{rav1}) and in the second word of (\ref{u58}).

In  (\ref{rav1}) with $k=n-3$
we have $s=n-2,$ $m=n;$
and  (\ref{rav1}) becomes
\begin{equation}
e_{n-3n}e_{n-3n-2}\equiv _{n-2}\beta e_{n-3n-1}
\underline{x_{n-3}x_{n-2}x_nx_{n-2}}.
\label{dav1}
\end{equation}
In view of Lemma  \ref{ong4d}, we may
transform the underlined part in  $U_n$
neglecting the words starting with  $x_{n-2}^2$
and $x_n$
in much the same way as in  (\ref{rav2}),
with  $x_n$ in place of $x_{k+1}$.
So (\ref{dav1}) reduces to the required form.

The second word of  (\ref{u58}) with  $k=n-3$
assumes the form
$e_{n-3n}^2x_{n-2}x_{n-1}^2=e_{n-3n}
\underline{x_{n-3}x_{n-2}x_nx_{n-2}}x_{n-1}^2.$
By Lemma \ref{ong4a} applied to  $U_n,$
the underlined word is a linear combination of words starting
with $x_{n-2}$
and $x_n.$
However, by Lemma  \ref{ong4d} both  $e_{n-3n}x_{n-2}$
and $e_{n-3n}x_n$
equal zero up to $\equiv _{n-2}.$ \hfill $\Box $
\addtocounter{nom}{1}
\begin{lemma}
The set $B$
satisfies the conditions of Lemma $\ref{tvp}.$
\label{dn7}
\end{lemma}
{\it Proof}. By Lemmas  \ref{dn6} and  \ref{kri1} one need show only
that in $U_P^b(\frak g)$
the words  (\ref{dGB}) are linear combinations of lesser ones.
The words  $v_6$ with $m=n-2,$ and $u_0,$ $u_1,$ $u_1^{\prime },$ $u_2$
have the required decomposition  since
they belong either to $U_{n-1}$
or to $U_{n}.$
Lemma  \ref{ong4d}
shows that $v_3\equiv _{k+1}0,$ $v_4\equiv _{k+1}0,$
$v_4^{\prime }\equiv _{k+1}0.$
Lemma \ref{d5} with $s=m-1$
yields the required representation for  $v_5.$
Consider $v_6$
with $m=n-1.$
Let us prove by downward induction on $k$
that
$$
u_{k\, n-1}e_{kn}\equiv _{k+1}\varepsilon e_{kn}u_{k\, n-1},\ \ \varepsilon \neq 0.
$$
For $k=n-1$
this equality assumes the form  (\ref{rD3}). Let  $k<n-1.$
Let us transpose the second letter $x_k$ of $u_{k\, n-1}e_{kn}$
as far to the left as possible in $U_{n-1}.$
After an application of (\ref{rA1}) we get
$$
u_{k\, n-1}e_{kn}\equiv _{k+1}
\alpha x_k^2(u_{k+1 n-1}e_{k+1n}), \ \ \ \alpha \neq 0.
$$
It suffices to  apply the inductive
supposition to the term in the parentheses,
and then by  (\ref{rA1}) and (\ref{rA3}) for  $U_n$
to move $x_k$
to the proper place. \hfill $\Box $
\addtocounter{nom}{1}
\begin{lemma}
If $p_{11}\neq 1$
then the values of  $[v]^h,$
where  $[v]\in B,$ $v\neq x_i,$ $h\geq 1$
are not skew primitive, in particular they are non-zero.
\label{dlpr}
\end{lemma}

{\it Proof}. One need consider only  super-letters
that belong neither to $U_{n-1}$
nor to $U_n.$
That is  $[e_{km}]$
with $m<n.$
We use induction on $n.$

For $n=3$
the algebra of the type  $D_3$
reduces to the algebra of the type  $A_3$
with a new ordering of variables $x_2>x_1>x_3.$
Therefore we may use Theorem {\bf $A_n,$}
after the decomposition below of $e_{12}$
in the PBW-basis:
$$[[x_1x_3]x_2]=-p_{12}p_{32}[x_2[x_1x_3]]
+\beta [x_1x_3]\cdot x_2.$$

Let $n>3.$
If $k>1$
then the inductive supposition works. For $k=1,$ $m>2$
we have $e_{1m}=[x_1[e_{2m}]],$
and one may  repeat the arguments
of Lemma  \ref{alpr} with $e$
in place of $u$
starting at (\ref{pov1}).
If $m=2$
then we may repeat the arguments of Lemma \ref{blpr} with $e$
on place of  $w$
starting  at (\ref{pov3}).  \hfill $\Box $

\

{\it Proof} of Theorem $D_n.$
For the first statement it will suffice to prove that
all super-letters  (\ref{DN1}) are hard in $U_P^b({\frak g}).$

Since none of  $u_{km}$
contains sub-words  (\ref{vor}), $[u_{km}]$ are hard.

Suppose  $[e_{km}]$
is non-hard. By Lemmas \ref{dn7} and  \ref{tvp}
all hard super-letters belong to $B.$
Thus, by Lemma \ref{mdn}, we get  $[e_{km}]=0.$
Since   deg$_n(e_{km})=$deg$_{n-1}(e_{km})=1,$
the equality  $[e_{km}]=0$
is also valid in the algebra $D^{\prime }$
defined by the same relations as
$U_P^b({\frak g})$
is, but $[x_{n-2}x_n^2]=0$ and  $[x_{n-2}x_{n-1}^2]=0.$
Let us equate to zero all monomials in all
the defining relations of  $D^{\prime },$
but  $[x_{n-1}x_n]=0.$
Consider the algebra  $R^{\prime }$
defined by  (\ref{rD3})  and by the resulting system
of monomial relations. It is easy to verify that
the mentioned relations system  $\Sigma $
of $R^{\prime }$
is closed under the compositions. Since $e_{km}$
contains none of leading words of $\Sigma ,$
the super-letter $[e_{km}]$
is non-zero in $R^{\prime },$
and so in $D^{\prime }$ too.
This contradiction proves the first statement.

If  $[v]^h, [v]\in B$
is of finite height then by Lemma \ref{mdn}
and the homogeneous version of Definition \ref{h1}
we have  $[v]^h=0.$
For $p_{11}\neq 1$
this contradicts Lemma  \ref{dlpr}.
In a similar manner, by Lemma  \ref{prim}, every skew
 primitive homogeneous  element has the form $[v]^h.$
This, together with Lemma  \ref{dlpr}, proves  both the fourth
statement and the second one with $p_{11}\neq 1.$

If $p_{11}=1$ then by  (\ref{suxd}) we have
$p_{ii}=p_{ij}p_{ji}=1$
for all $i,j.$
This means that the skew commutator itself is
a quantum operation. Hence all elements of  $B$
are skew-primitive. These elements span a colour Lie super-algebra.
Now, as in Theorem  $A_n,$
one may use the PBW theorem for colour Lie super-algebras.

For the third statement it will suffice to show that all
super-letters (\ref{dGB}) are zero in  $U_P({\frak g}).$
We have proved already that they are non-hard.
Therefore it remains to use the homogeneous version of Definition \ref{tv1}
and Lemma \ref{mdn}.     \hfill $\Box $

\

\centerline {8. \it Conclusion}

\

\stepcounter{par}
\noindent
We see that in all Theorems $A_n$--$D_n$
the lists of hard super-letters are independent of the parameters $p_{ij}.$
This fact signifies that the Lalonde--Ram basis of
the ground Lie algebra (see, [{\bf \ref{pLR}}, Figure 1]) with the skew
commutator in place of the Lie operation coincides
with the set of all hard super-letters.
It is very interesting to clarify how general this statement is.
On the one hand, this does
not hold without exception for all quantum enveloping
algebras since in Theorems
$A_n$--$D_n$ a restriction does exist. If $p_{ii}=-1,$ $1\leq i<n,$ $n>2$ then it is easy
to see by means of the Diamond Lemma that the sets of hard super-letters are infinite. On
the other hand, this is not a specific
property of Lie algebras defined by the Serre
relations. By the
Shirshov theorem [{\bf \ref{pSh1}}]
any relation can be reduced to a linear combination of
standard nonassociative words.
\addtocounter{nom}{1}
\begin{corollary}
If ${\frak g}$ is defined by the only relation $f=0,$
where $f$ is a linear combination of
standard nonassociative words,
then the set of all hard in $U_P({\frak g})$ super-letters
coincides with the Hall--Shirshov basis
of  ${\frak g}$  with the skew commutator in place of the Lie operation.
\label{fi}
\end{corollary}
{\it Proof}.
The only  relation $f^*=0$ forms a Groebner--Shirshov  system
since, according to 1s,
none of onsets of its leading word, say $w,$
coincides with a proper terminal of $w$. Consequently,
a super-letter $[u]$ is hard if and
only if $u$ does not contain $w$ as a sub-word.
We see that this criteria is independent
of $p_{ij}$ as well. \hfill $\Box $

\

Furthermore, the third statement of Theorem $A_n$
shows that $U_P^b({\frak g})$
can be defined by the following relations
in the PBW-generators $X_u=[u].$
\begin{equation}
\begin{matrix}
[X_u,X_v]=0, \hfill & u>v,\ &[[u][v]]\notin B\hfill \cr
[X_u,X_v]=X_{uv}, \hfill &  &[[u][v]]\in B.\hfill \end{matrix}
\label{Ank}
\end{equation}
This is an argument in favour of considering
the super-letters PBW-generators
{\bf k}$[G]$-module as a quantum analogue of a Lie algebra. However in the cases
$B_n,$ $C_n,$ $D_n$ the defining relations became more complicated.
For example,
\begin{equation}
\begin{matrix} B_n:\ \ [[u_{k\, n-1}][w_{kn}]]=\alpha [u_{kn}]^2,\hfill & \alpha \neq 0 \hbox{ if
} p_{nn}\neq 1;\hfill \cr
C_n:\ \ [[u_{k\, n-2}][v_{k\, n-1}]]=\alpha [v_k]+\beta [u_{kn}]\cdot [u_{k\, n-1}], \hfill &
\beta \neq 0 \hbox{ if } p_{11}\neq 1;\hfill \cr
D_n:\ \ [[u_{k\, n-2}][e_{k\, n-1}]]=\alpha [e_{kn}]\cdot [u_{k\, n-1}],\hfill & \alpha \neq 0
\hbox{ if } p_{11}\neq \pm 1.\hfill \end{matrix}
\label{comp}
\end{equation}

It is far more interesting that for $p_{11}\neq 1$
the algebra ${\frak g}_P$
turns out to be very simple in structure.
Only unary quantum operations can be non-zero.
Other ones may be defined, but due to the homogeneity
their values equal zero. In
particular, if $p_{11}^{[t]}\neq 0$
then without exception all quantum operations have zero values. This provides reason
enough to consider $U_P({\frak g})=$ $U({\frak g}_P)$
as an algebra of  `commutative' quantum polynomials.
Certainly it is very interesting to
elucidate to what extent this statement is still
retained for  the quantum universal
enveloping algebras of homogeneous
components of other Kac--Moody algebras defined by the
Gabber--Kac relations (\ref{rel}). Also it is interesting
to investigate the structure of other
`commutative' quantum polynomial algebras. For example,
one may note that if a semi-group generated by $p_{ij}p_{ji}$
does not contain 1, then
$G\langle x_1,\ldots ,x_n \rangle $
itself is a `commutative' quantum polynomial algebra merely since in this case there
exists no non-zero quantum operation at all. In another extreme case when
$p_{ij}p_{ji}=1$ for all $i,j,$ the `commutative' quantum variables commute by
$x_ix_j=p_{ij}x_jx_i.$

In a similar manner, the Drinfeld--Jimbo enveloping algebra can be considered as a
`quantum' Weyl algebra of (skew) differential operators (see Sec. 6). The resulting
`quantum' Weyl algebra is simple in the following sense.
\addtocounter{nom}{1}
\begin{corollary}
Let ${\frak g}$ be a simple finite dimensional
Lie algebra of the infinite series. If
$q^{[m]}\neq 0,$ $m\geq 2$ then every non-zero Hopf ideal $I$
of the Drinfeld--Jimbo enveloping
algebra contains all generators $x_i, x^-_i.$
\label{fin}
\end{corollary}
{\it Proof}. By the Heyneman--Radford theorem,
the ideal $I$ has a non-zero skew primitive
element, say $a.$ According to Lemma \ref{raz} and Theorems  $A_n$--$C_n,$
the element $a$ is either a constant, $\alpha (1-g),$ or proportional to one of the elements
$x_i, x^-_i.$ In the former case $I$ contains all $x_i$ with $\chi  ^i(g)\neq 1$ since
$x_i a-\chi ^i(g)a_ix_i=\alpha (1-\chi ^i(g))x_i.$
Here the equality $\chi ^i(g)=1$ can not be valid for all $i$ since
$\chi ^i(g_j)=q^{-d_ia_{ij}}$ (see, Example 4 of Section 2) and the columns of the Cartan
matrix are linearly independent. In the latter case (and now in the former one as well) we
get $[x_i,x_i^-]=\varepsilon _i(1-g_i^2)\in I,$
i.e. as above $I$ contains all elements $y=x_i^{\pm }$  with $1\neq \chi ^y(g_i^2)$
$=q^{\pm 2d_ja_{ij}}.$ Since the Coxeter graph is connected, $I$ contains all
$x_i, x^-_i.$     \hfill $\Box $

\

{\it Acknowledgements.} The author is grateful to Dr. J. A. Montaraz,
the director of the FES-C UNAM,  Dra. S. Rodr\'\i guez-Romo, and A.V. Lara Sagahon
for providing  facilities  for the research and also to Dr. L.A. Bokut'
and Dr. R. Bautista for helpful comments  on the subject matter.

\

\centerline{\large \it References}

\

\newcounter{items}
\begin{list}{\bf \arabic{items}.}{\usecounter{items}}

  \item  \label {pBau} {\sc C. Bautista,} `A Poinkare--Birkhoff--Witt
theorem for generalized Lie color algebras', {\it Journal of Mathematical
Physics} 39 N7(1998) 3829--3843.

  \item  \label {pBMM} {\sc K.I. Beidar, W.S. Martindale III,}
and {\sc A.V. Mikhalev,}
{\it Rings with Generalized Identities} (Pure and Applies Mathematics
196, Marcel Dekker, New York--Basel--Hong Kong, 1996).

  \item  \label {pBer} {\sc G.M. Bergman,} `The diamond lemma for ring theory',
{\it Adv. in Math.} 29 N2(1978) 178--218.

 \item  \label {pBo} {\sc L.A. Bokut',}
`Unsolvability of the word problem and subalgebras of finitely presented Lie algebras',
{\it Izv.Akad.Nauk. Ser. Mat.} 36 N6(1972) 1173--1219.

  \item  \label {pBok} {\sc L.A. Bokut',} `Imbeddings into simple
associative algebras',
{\it Algebra and Logic} 15 N2(1976) 117--142.

  \item  \label {pKle} {\sc L.A. Bokut',} and {\sc A.A. Klein,}
`Serre relations
and Groebner--Shirshov bases for simple Lie algebras I, II',
{\it International Journal of Algebra and Computation}
 6 N4(1996) 389--412.

 \item  \label {pBK} {\sc L.A. Bokut',} and {\sc G.P. Kukin,}
{\it Algoritmic and Combinatorial Algebra} (Mathematics and Its Applications
255, Kluwer Academic Publishers, Dordrecht-Boston-London, 1994).

  \item  \label {pBM} {\sc L.A. Bokut',} and {\sc P. Malcolmson,} `Groebner bases for quantum
enveloping algebras', {\it Israel Journal of Mathematics}
96(1996) 97--113.

 \item  \label {pBor} {\sc R. Borcherds,} `Generalized Kac-Moody algebras', {\it Journal of Algebra}, 11(1988) 501--512.

  \item  \label {pFox} {\sc K.T. Chen, R.H. Fox,} and {\sc R.C. Lyndon,}
`Free differential
calculus IV, the quotient groups of the lower central series',
{\it Ann. of Math.} 68(1958) 81--95.

  \item  \label {pCli} {\sc G. Clift,} `Crystal bases and Young tableaux',
{\it Journal of Algebra} 202 N1(1998) 10--35.

  \item  \label {pCohn} {\sc P.M. Cohn,} `Sur le crit\`ere de Friedrichs
pour les commutateur dans une alg\`ebre associative libre',
{\it C. r. Acad. sci. Paris} 239 N13(1954) 743--745.

  \item  \label {pCoh} {\sc P.M. Cohn,} {\it Universal Algebra}
(Harper and Row, New-York, 1965).

  \item  \label {pDri} {\sc V.G. Drinfeld,} `Hopf algebras and the Yang--Baxter
equation', {\it Soviet Math. Dokl.}, 32(1985) 254--258.

  \item  \label {pFri}  {\sc K.O. Friedrichs,} `Mathematical aspects
of the quantum theory of fields. V', {\it Communications in Pure and
Applied Mathematics} 6(1953) 1--72.

  \item  \label {pKac} {\sc O.Gabber,} and {\sc V. Kac,} `On defining relations of  certain
infinite-dimensional Lie algebras', {\it Bulletin (New series)
of the American  Mathematical Society} 5, N2(1981) 185--189.

  \item  \label {pJim} {\sc M. Jimbo,} `A q-difference analogue
of $U(\frak g)$
and the Yang--Baxter equation', {\it Lett. Math. Phis.}
10(1985) 63--69.

 \item  \label {pKang} {\sc S.-J. Kang,} `Quantum deformations of generalized Kac-Moody algebras and their modules',
{\it Journal of Algebra}, 175(1995) 1041--1066.

  \item  \label {pKas}  {\sc M. Kashiwara,} `Crystallizing the $q$-analogue
of universal enveloping algebras', {\it Comm. Math. Phis.}
133(1990) 249--260.

  \item  \label {pKas1}  {\sc M. Kashiwara,} `On crystal bases of the q-analog of
universal enveloping algebras', {\it Duke Mathematical Journal}
63 N2(1991) 465--516.

  \item  \label {pKh} {\sc V.K. Kharchenko,} `An algebra of skew
primitive elements', {\it Algebra and Logic} 37 N2(1998) 101--126.

  \item  \label {pKh1} {\sc V.K. Kharchenko,}
`A quantum analogue of the Poincar\`e--Birkhoff--Witt theorem',
{\it Algebra and Logic} 38 N4(1999) 476--507; English translation 259--276.

 \item  \label {pJA} {\sc V.K. Kharchenko,}
` An existence condition for multilinear quantum operations',
{\it Journal of Algebra} 217(1999) 188--228.

 \item  \label {pDAN} {\sc V.K. Kharchenko,}
` Character Hopf algebras and quantizations of Lie algebras',
{\it Doklady Mathematics} 60 N3(1999) 328--329.

  \item  \label {pKit} {\sc A. Kuniba, K.C. Misra, M. Okado, T. Takagi,}
and {\sc J. Uchiyama,} `Crystals for Demazure modules of classical
affine Lie algebras', {\it Journal of Algebra} 208(1998) 185--215.

  \item  \label {pLR} {\sc M. Lalonde,} and {\sc A. Ram,} `Standard Lyndon
bases of Lie algebras and enveloping algebras',
{\it Trans. Amer. Math. Soc.} 347 N5(1995) 1821--1830.

  \item  \label {pLot} {\sc M. Lothaire,} {\it Combinatorics on words},
(Encyclopedia of Mathematics and its Applications 17,
Addison--Wesley Publ. Co. 1983).

  \item  \label {pLu} {\sc G. Lusztig,} `Quantum groups at roots of 1',
{\it Geometria Dedicada} 35, N1-3(1990) 89--113.
  \item  \label {pLuz} {\sc G. Lusztig,}  {\it Introduction to Quantum Groups}
(Progress in Mathematics 10, Birkhauser Boston, 1993).

  \item  \label {pLyn} {\sc R.C. Lyndon,} `A theorem of Friedrichs',
{\it Michigan Mathematical Journal} 3, N1(1955--1956) 27--29.

  \item  \label {pLyu} {\sc V. Lyubashenko,} and {\sc A. Sudbery,} `Generalized Lie algebras
of type $A_n$', {\it Journal of Mathematical Physics} 39,
N6(1998) 3487--3504.

  \item  \label {pMag} {\sc W. Magnus,} `On the exponential solution of differential
equations for a linear operator', {\it Communications in Pure and
Applied Mathematics}  7(1954) 649--673.

  \item  \label {pMMo} {\sc J.W. Milnor} and {\sc J.C. Moore,}
`On the structure
of Hopf algebras', {\it Annals of Math.} 81(1965) 211--264.

  \item  \label {pMon} {\sc S. Montgomery,} {\it Hopf Algebras and Their Actions
on Rings} (CBMS 82, AMS, Providence, 1993).

  \item  \label {pRad} {\sc D.E. Radford,} `The structure of Hopf algebras
with projection',
{\it Journal of Algebra} 92(1985) 322--347.

  \item  \label {pSh1} {\sc A.I. Shirshov,} `On free Lie rings',
{\it Matem. Sbornic} 45(87) N2(1958) 113--122.

  \item  \label {pSh2} {\sc A.I. Shirshov,} `Some algorithmic problems
for Lie algebras',
{\it Sibirskii Math. Journal} 3 N2(1962) 292--296.

  \item  \label {pYa} {\sc Yamane,} `A Poincar\`e-Birkhoff-Witt theorem
for quantized universal enveloping algebras of type $A_N$',
{\it Publ. RIMS. Kyoto Univ.} 25(1989) 503--520.
\end{list}

\end{document}